\documentclass[11pt,leqno]{amsart}

\usepackage{amsmath}
\usepackage{amsthm}
\usepackage{graphicx}
\usepackage{mathtools}
\usepackage{amsfonts}
\usepackage{amssymb}
\usepackage{hyperref}
\usepackage{bm}
\usepackage[margin=1 in]{geometry}
\usepackage{enumerate}
\usepackage[normalem]{ulem}
\usepackage{tikz}
\usepackage{tikz-cd}
\usepackage{xcolor}
\usetikzlibrary{matrix,arrows,positioning,automata}
  \usepackage{cancel}
  \usepackage{todonotes}
\usepackage{cite}
  
\numberwithin{equation}{section}

\newcommand{\R}{\mathbb R}
\def\E{\mathbb E}

\def\XXint#1#2#3{{\setbox0=\hbox{$#1{#2#3}{\int}$}
\vcenter{\hbox{$#2#3$}}\kern-.5\wd0}}

\newcommand{\T}{\mathbb{T}}

\numberwithin{equation}{section}
\newtheorem{thm}{Theorem}[section]
\newtheorem{lem}[thm]{Lemma}
\newtheorem{cor}[thm]{Corollary}
\newtheorem{prop}[thm]{Proposition}

\theoremstyle{definition}
\newtheorem{defn}[thm]{Definition}
\newtheorem{rmk}[thm]{Remark}
\newtheorem{ex}[thm]{Example}

\def\smallnegint{\mathop{\int\mkern-13mu
        \raise.5ex\hbox{${\scriptscriptstyle\diagup}$}}\nolimits}

\def\div{\operatorname{div}}
\def\tr{\operatorname{tr}}

\def\ssetminus{\,\raise.4ex\hbox{$\scriptstyle\setminus$}\,}

\newcommand{\be}{\begin{equation}}
\newcommand{\ee}{\end{equation}}

\newcommand{\bc}{\begin{case}}
\newcommand{\ec}{\end{cases}}
\newcommand{\bs}{\begin{split}}
\newcommand{\es}{\end{split}}

\newcommand{\norm}[1]{\left\Vert#1\right\Vert}

\renewcommand{\d}{d}

\newcommand{\Rd}{{\mathbb{R}^\d}}

\renewcommand{\bar}{\overline}
\renewcommand{\tilde}{\widetilde}

\newcommand{\RR}{\mathbb{R}}

\newcommand{\ac}{\mathrm{ac}}
\newcommand{\del}{\partial}
\newcommand{\Id}{\mathrm{Id}}
\newcommand{\Lip}{\mathrm{Lip}}
\newcommand{\loc}{\mathrm{loc}}

\newcommand{\mcl}{\mathcal}
\newcommand{\oline}{\overline}
\newcommand{\oo}{\infty}
\newcommand{\Tan}{\mathrm{Tan}}

\title[A comparison principle for HJB equations in the Wasserstein space]{A comparison principle for semilinear Hamilton-Jacobi-Bellman equations in the Wasserstein space}

\date{\today}

\begin{document}

\author[S. Daudin]{Samuel Daudin 
\address{(S. Daudin) Universit\'e C\^ote d'Azur, CNRS, Laboratoire J.A. Dieudonn\'e, 06108 Nice, France
}\email{samuel.daudin@univ-cotedazur.fr
}}
 
\author[B. Seeger]{ Benjamin Seeger
\address{(B. Seeger) University of Texas at Austin, 2515 Speedway, PMA 8.100, Austin, TX 78712, USA
}\email{seeger@math.utexas.edu
}}

\thanks{S. Daudin acknowledges the financial support of the European Research Council (ERC) under the European Union's Horizon 2020 research and innovation program (ELISA project, Grant agreement No. 101054746). Part of this research was achieved while he was affiliated to the CEREMADE, Universit\'e Paris-Dauphine. B. Seeger was partially supported by the National Science Foundation (NSF) under award number DMS-1840314. Both authors acknowledge support from the IMSI (University of Chicago) where they were invited during the long program Distributed Solutions to Complex Societal Problems.}

\maketitle

\begin{abstract}
    The goal of this paper is to prove a comparison principle for viscosity solutions of semilinear Hamilton-Jacobi equations in the space of probability measures. The method involves leveraging differentiability properties of the $2$-Wasserstein distance in the doubling of variables argument, which is done by introducing a further entropy penalization that ensures that the relevant optima are achieved at positive, Lipschitz continuous densities with finite Fischer information. This allows to prove uniqueness and stability of viscosity solutions in the class of bounded Lipschitz continuous (with respect to the $1$-Wasserstein distance) functions. The result does not appeal to a mean field control formulation of the equation, and, as such, applies to equations with nonconvex Hamiltonians and measure-dependent volatility. For convex Hamiltonians that derive from a potential, we prove that the value function associated with a suitable mean-field optimal control problem with nondegenerate idiosyncratic noise is indeed the unique viscosity solution.
\end{abstract}

\tableofcontents




\section{Introduction}

This paper proposes a novel approach to the well-posedness of semilinear Hamilton-Jacobi-Bellman equations on the space of probability measures, taking the form
\begin{equation}\label{eq:HJB}
	\left\{
	\begin{split}
		&-\del_t U(t,\mu) - \int_{\Rd} \tr\left[a(\mu,x) \nabla_x \partial_\mu U(t,\mu,x)\right]d\mu(x) \\
		& \qquad + \int_{\Rd} H(x, \partial_\mu U(t,\mu,x), \mu)d\mu(x) = 0 \quad \text{in } [0,T] \times \mcl P_2(\R^d), \\
	&U(T, \mu) = \mathcal{G}(\mu) \quad \text{in } \mcl P_2(\R^d),
	 \end{split}
	\right.
\end{equation}
where $\mathcal{P}_2(\R^d)$ is the Wasserstein space of probability measures over $\R^d$ with finite second order moment and $T>0$ is a finite horizon. The unknown $U$ is defined over $[0,T] \times \mathcal{P}_2(\R^d)$, $\partial_{\mu}U(t,\mu,.) \in L^2(\mu; \R^d)$ is the Wasserstein gradient of $U(t,.)$ at $\mu$ and $\nabla_x \partial_\mu U$ is its gradient in the additional variable $x$. The Hamiltonian $H$ is defined over $\R^d \times \R^d \times \mathcal{P}_2(\R^d)$, the diffusion matrix $a(\mu,x)$ is positive definite for all $(\mu,x) \in \mathcal{P}_2(\R^d) \times \R^d$, and $\mathcal{G} : \mathcal{P}_2(\R^d) \rightarrow \R$ is the terminal data (precise definitions and assumptions will be given later on).

Equations of the form \eqref{eq:HJB} have recently attracted much attention due to their connection with mean-field stochastic optimal control and potential mean-field games. Indeed, when the Hamiltonian $H$ is convex in the $p$-variable and takes the form $H(x,p,\mu) = \sup_{q \in \R^d} \bigl \{ -p \cdot q - L(x,q,\mu) \bigr \}$ for some Lagrangian $L : \R^d \times \R^d \times \mathcal{P}_2(\R^d)$, the unique solution to \eqref{eq:HJB} is expected to be given, for every $(t_0,\mu_0) \in [0,T] \times \mathcal{P}_2(\R^d)$ , by
\begin{equation} 
 U(t_0,\mu_0) = \inf_{(\mu,\alpha)} \int_{t_0}^T \int_{\R^d} L \bigl(x,\alpha_t(x), \mu_t \bigr)d\mu_t(x)dt + \mathcal{G}(\mu_T) 
 \label{valuefunctionintro2023}
\end{equation}
where the infimum runs over couples $(\mu,\alpha)$ such that $\mu \in \mathcal{C}([t_0,T] ,\mathcal{P}_2(\R^d))$, $\alpha \in L^2([0,T] \times \Rd, dt \otimes \mu_t; \Rd)$, and the following (non-local) Fokker-Planck equation, describing the law of the solution to a McKean-Vlasov stochastic differential equation, is satisfied:
\begin{equation}
\left \{
\begin{array}{ll}
\partial_t \mu + \div((\alpha_t(x)\mu) -\partial^2_{ij}\bigl( a^{ij}(x,\mu_t) \mu \bigr) = 0 &\mbox{ in } (t_0,T) \times \R^d,  \\
\mu(t_0)=\mu_0.
\end{array}
\right.
\label{eq:FPEintro2023}
\end{equation}

Even if the diffusion matrix $a$ is non-degenerate and the data is smooth, the value function $U$ is not necessarily differentiable (see Section \ref{sec:exmaples} below); equivalently, for a given $(t_0,\mu_0) \in [0,T] \times \mcl P_2(\Rd)$, different optimal controls can be found for the associated control problem. For this reason, it is necessary to introduce a notion of weak solutions to \eqref{eq:HJB}.

In this paper we propose a direct approach to the comparison principle for viscosity sub and supersolutions of \eqref{eq:HJB}, using the doubling of variables argument of the classical theory of viscosity solutions, in which the squared $2$-Wasserstein distance takes the role of a penalizing test function. To ensure that relevant optima are achieved at sufficiently regular measures where the $2$-Wasserstein distance is differentiable and can be used as a test function in the equation, we introduce a further entropy penalization. This allows us to prove comparison and stability in the class of bounded functions that are Lipschitz continuous with respect to the $1$-Wasserstein distance. Notably, our method for uniqueness does not rely on the control formulation, and therefore applies to non-convex Hamiltonians. When the Hamiltonian derives from a Lagrangian, we show that the value function \eqref{valuefunctionintro2023} is indeed the unique viscosity solution.

One of the key technical steps of the proof consists in showing regularity properties of maxima of expressions of the form
\begin{equation} 
\mathcal{P}_2(\R^d) \times \mathcal{P}_2(\R^d) \ni (\mu,\nu) \mapsto U(\mu) - V(\nu) - \frac{1}{2\epsilon}d_2^2(\mu,\nu) - \delta \mathcal{E}(\mu) - \delta \mathcal{E}(\nu)
\label{MaximumdoublingIntro2023}
\end{equation}
where $U$ and $V$ play respectively the roles of sub and super-solutions, $d_2$ is the $2$-Wasserstein distance, $\mathcal{E}(\sigma) = \int_{\R^d} \log \bigl( \sigma(x) \bigr) \sigma(x)dx $ is the entropy (set to be $+\infty$
if $\sigma$ is not absolutely continuous with respect to the Lebesgue measures) of $\sigma \in \mathcal{P}_2(\R^d)$ and $\epsilon,\delta >0$ are small parameters. In particular, the further entropy penalizations force the optima $(\mu,\nu)$ to be attained at sufficiently smooth densities with full support, and whose unique transport map is a diffeomorphism. To show this, we rely on techniques similar to those developed in the analysis of the celebrated Jordan-Kinderlehrer-Otto scheme for the Fokker-Planck equation, see \cite{Jordan1998,Sant_15}. In this context the discrete scheme involves solving variational problems of the form 
\begin{equation} 
\max \Bigl \{ \mathcal{P}_2(\R^d) \ni \mu \mapsto  \int_{\R^d}u(x)d\mu(x) - \frac{1}{2\epsilon} d_2^2(\mu,\nu) - \delta \mathcal{E}(\mu)  \Bigr \}, 
\end{equation}
where $\nu$ is a given measure and $u: \R^d \rightarrow \R$ is a smooth function. One challenge we face is to find minimal regularity conditions on $U$ and $V$ that still ensure enough regularity for the maxima of \eqref{MaximumdoublingIntro2023}, and that are satisfied by the candidate solution \eqref{valuefunctionintro2023}. We find that requiring $U$ and $V$ to be Lipschitz continuous with respect to the $1$-Wasserstein distance is enough for our purpose. We refer to Section \ref{sec:NotionofSolandMainresults} for a detailed discussion on the strategy of the proof. Throughout the paper, we invoke the differentiability properties of the $2$-Wasserstein distance that were put forward in \cite{Gigli2004}, see also \cite{AGS_2008, Alfonsi2020}, as well as on regularity results for the Brenier transport map, first obtained in \cite{Caffarelli1991} and further extended in \cite{Cordero2019}. The use of entropy penalization in this spirit has also been used in works of Feng and co-authors to study large deviations \cite{FK_book_06} and to prove comparison for a different class of Hamilton-Jacobi equations deriving from action integrals related to certain systems of conservation laws \cite{Feng2012}.

Mean-field control problems like \eqref{valuefunctionintro2023} arise formally as limits of control problems involving a large number of controlled particles subject to independent idiosyncratic noises. The convergence problem consists of verifying this limit rigorously and further understanding further qualitative and quantitative information, such as the rate of convergence, propagation of chaos for the optimal trajectories, and convergence of the controls; see \cite{Lacker2017, CDJS_23, Cardaliaguet2023, DDJ_23}. The recent contributions \cite{CDJS_23, Cardaliaguet2023, DDJ_23} show that regularity of the value function \eqref{valuefunctionintro2023}, as well as its satisfying the viscosity solution property, are central to the understanding of the convergence problem. 

Another approach to analyzing nonlinear equations on Wasserstein spaces is to “lift'' the equation to a Hilbert space of random variables, a strategy which was introduced by Lions in his lectures at the Coll\`ege de France \cite{PLLCdF} to study classical solutions of master equations from the theory of mean field games. This viewpoint has since been used to study certain variants of the HJB equation \eqref{eq:HJB} \cite{GT_19, MS_23} (see also \cite{B_23} for a recent related approach), or for first- and second-order mean field game master equations with common noise and no idiosyncratic noise \cite{CS_22}, by adapting the theory of Crandall and Lions \cite{Crandall1992,PLL1988,PLL1989,CL_85,CL_86} for viscosity solutions to fully non-linear equations on infinite-dimensional spaces. A major difference with our setting is the nature of the second-order term in \eqref{eq:HJB}, which arises from the mean field control problem with idiosyncratic noise, rather than common noise (or zero noise), as those studied in the aforementioned works. Certain technical obstructions seem to place the Hilbertian method out of reach for \eqref{eq:HJB} and other equations involving idiosyncratic noise (see Remark \ref{rmk:HilbertianApproach} below for a detailed discussion), and, therefore, much effort has been made to prove well-posedness for an “intrinsic'' notion of viscosity solutions that does not appeal to the lifting procedure. 

The main result in this direction can be found in \cite{Cosso2021}, in which existence and uniqueness was proved for viscosity solutions of a certain class of fully non-linear equations arising from mean-field control problems with controlled, possibly degenerate (but not measure-dependent) volatility. The general strategy of \cite{Cosso2021} is to compare sub and supersolutions to certain smooth approximations of the value function, and, in particular, it relies on the finite-dimensional counterpart of the control problem \eqref{valuefunctionintro2023} with a large, but finite, number of interacting particles. Another approach developed recently in \cite{Soner2022} consists in replacing the $2$-Wasserstein distance in the doubling of variables argument with weaker, negative Sobolev-type norms that enjoy better differentiability properties. This allows for the comparison of viscosity sub/super solutions, as soon one of them is Lipschitz continuous with respect to this Sobolev topology. Although some regularity is assumed for the sub and supersolutions in our proof of the comparison principle, the requirement (Lipschitz continuity in the $d_1$-distance) is much weaker than that used in \cite{Soner2022}. On the other hand, our approach is limited to semilinear equations with non-degenerate diffusions.

As we have already mentioned, our approach, especially compared with \cite{Cosso2021}, does not rely on the control formulation, and, as such, allows for the treatment of equations with general, non-convex Hamiltonians. This raises the potential for the use of our method for Isaacs-type equations arising from differential games, see \cite{Cosso2019}. We also expect that our method can accommodate equations corresponding to control problems for which the $N$-particle problem does not give smooth approximations of the solution to the infinite-dimensional mean-field control problem. This is for instance the case when the controlled equation \eqref{eq:FPEintro2023} faces additional state constraints, see \cite{Daudin2023,Daudin2021}, and we expect that the value function in that setting is the unique constrained solution to the HJB equation. We aim to address these settings in future work.

Let us also mention two other approaches to the well-posedness of equation \eqref{eq:HJB}. In \cite{Cecchin2022}, the authors prove well-posedness in the class of Lipschitz (with respect to the metric of $\bigl(\mathcal{C}^2 (\T^d) \bigr)^*$) and (displacement) semi-concave functions after reformulating the equation as an equation over the Fourier coefficients of the probability measures. The equation is then understood in a suitable almost-everywhere sense, similar to early works on finite dimensional HJB equations. This formulation is then sufficient to give a conservative meaning to the Master equation from mean-field game theory which is satisfied by the linear functional derivative of \eqref{valuefunctionintro2023}. The strategy and the difficulties in \cite{Cecchin2022} are, in some sense, orthogonal to those in the present paper. Indeed, formulating \eqref{eq:HJB} over the space of Fourier coefficients allows for a simple way to deal with the second order linear term, and most of the difficulty comes from the nonlinear Hamiltonian term. In both \cite{Cecchin2022} and our setting, the challenge is to reduce the analysis of \eqref{valuefunctionintro2023} to interior points (in a suitable sense) of the space of probability measures. In \cite{Conforti2021} the authors obtain uniqueness (but not existence) for solutions to general equations in metric spaces that arise as gradient flows of functionals verifying a certain evolutionary variational inequality. For particular choices of $H$ and $\mathcal{G}$, equation \eqref{eq:HJB} can be recast in this framework.

Once \eqref{eq:HJB} is known to be well-posed, a natural question is whether the solution enjoys further regularity properties. Such properties were thoroughly investigated in \cite{CDJS_23, DDJ_23, Cardaliaguet2023} for the value function \eqref{valuefunctionintro2023}, as a starting point toward quantitative  rates for the convergence problem in mean-field control. These results show that equation \eqref{eq:HJB} (at least in the more restrictive setting of our existence Theorem \ref{thm:ExistenceThm}) exhibits similar regularity properties to viscosity solutions to first order finite dimensional HJB equations. In particular, the value function is globally Lipschitz and semi-concave, and, moreover, the results of \cite{Cardaliaguet2023} show that the value function is differentiable in an open and dense set of $[0,T] \times \mathcal{P}_2(\R^d)$ and, in this set, the equation is satisfied in the strong sense. These results rely on the representation formula \eqref{valuefunctionintro2023} for the solution as well as the corresponding set of optimality conditions, which takes the form of a mean-field game system of PDEs (see the aforementioned references). One potential application of our direct proof of the comparison principle is the possibility to prove such regularity properties directly from the equation. In this way, we demonstrate in Section \ref{sec:UniquenessSection} how such a proof leads to a posteriori Lipschitz regularity (in the $d_2$-metric) for the unique solution. 


A related question is the potential regularizing effect of the non-degenerate diffusive term. Although no global gain of differentiability with respect to the measure argument can be expected (see Section \ref{sec:exmaples}), it is proved in \cite{Cardaliaguet2023} that, whenever the value function \eqref{valuefunctionintro2023} is differentiable with respect to the measure variable, then the Wasserstein gradient enjoys further regularity in the additional variable. Another consequence of the non-degeneracy, in the setting of Theorem \ref{thm:ExistenceThm}, is that, if the terminal condition $\mcl G$ is $d_1$-Lipschitz, then the value function is Lipschitz with respect to weaker negative-Sobolev type norms at any earlier time \cite[Proposition 3.4]{DDJ_23}. We refer to \cite{CST2022} for similar questions regarding linear equations in the space of probability measures.

The rest of the paper is organized as follows. In Section \ref{sec:NotionofSolandMainresults} we give our notion of viscosity solutions and our main results. We also recall some well-known results on the optimal transport problem for $d_2$ and present the strategy of the proof. In Section \ref{sec:EntropyRegularization} we study the variational problem \eqref{MaximumdoublingIntro2023} which appears in the doubling of variables argument. In Section \ref{sec:UniquenessSection} we give the proof of the comparison principle. We also obtain stability results as well as Lipschitz estimates for the unique solution. Finally in Section \ref{sec:ExistenceSection} we prove our main existence result.

\subsection*{Notation}

Let $\mathcal{P}(\R^d)$ be the set of Borel probability measures over $\R^d$. For any $p  \in [1,\oo)$, we let
\[
\mcl P_p(\Rd) = \left\{ \mu \in \mcl P(\R^d) : \mcl M_p(\mu) := \int_\Rd |x|^p d\mu(x) < \infty \right\},
\]
which becomes a complete metric space when equipped with the $p$-Wasserstein distance defined, for $\mu,\nu \in \mcl P_p(\R^d)$,
\[
	d_p(\mu,\nu) := \left(\inf_{\gamma \in \Gamma(\mu,\nu)} \iint_{\Rd \times \Rd} |x-y|^p d\gamma(x,y) \right)^{1/p},
\]
where $\Gamma(\mu,\nu) \subset \mcl P(\Rd \times \Rd)$ denotes the set of measures $\gamma$ with marginals $\mu$ and $\nu$. When a functional $F: \mcl P_p(\Rd) \to \RR$ is Lipschitz continuous with respect the distance $d_p$, we denote
\[
	\Lip(F;d_p) := \sup_{\mu,\nu \in \mcl P_p(\Rd), \; \mu \ne \nu } \frac{ |F(\mu) - F(\nu)|}{d_p(\mu,\nu)}.
\]

We also define
\[
	\mcl P_{p,\ac}(\R^d) := \left\{ \mu \in \mcl P_p(\R^d) : d\mu(x) = f(x)dx \text{ for some } f \in L^1(\R^d)\right\},
\]
and we often abuse notation and denote in a same way an element of $\mathcal{P}_{p,ac}(\R^d)$ and its density with respect to the Lebesgue measure.

For $p \in [1,+\infty]$, a measure space $(X,\mu)$, and a normed space $V$, we denote the following subspace of $\mu$-measurable functions $f: X \to V$ by
\[
	L^p(X,\mu; V) := \left\{ f : \norm{f}_{L^p(X,\mu; V)} := \left[ \int_X \norm{f(x)}_V^p d\mu(x) \right]^{1/p} < \oo \right\},
\]
with the obvious modification if $p = +\oo$. We often omit the range $V$ when it is clear from context, and, for $X = \Rd$, we usually write $L^p(\RR^d, \mu;V) = L^p(\mu;V)$. In the special case where $X = [0,T] \times \RR^d$ and the measure $\mu$ is given by
\[
	\int_{[0,T] \times \Rd} f(t,x) d\mu(t,x) = \int_0^T \left[ \int_\Rd f(t,x) d\mu_t(x) \right] dt, \quad \forall f \in \mathcal{C}_b \bigl([0,T] \times \R^d \bigr), 
\]
for some curve $[0,T] \ni t \mapsto \mu_t \in \mcl P(\Rd)$, we use the notation $L^p(X,\mu) = L^p([0,T] \times \RR^d , dt \otimes \mu_t)$.

For $n = 1,2,\ldots$, $p \in [1,+\oo]$, and an open set $U \subset \Rd$, we denote by $W^{n,p}(U)$ the usual Sobolev space of functions on $U$ with all distributional derivatives through order $n$ belonging to $L^p(U)$. The space $W^{n,p}_\loc = W^{n,p}_\loc(\Rd)$ consists of functions $f$ on $\Rd$ such that $f|_{U} \in W^{n,p}(U)$ for all bounded $U \subset \Rd$. Finally, we use the notation $H^n = W^{n,2}$.

\section{Notion of solution and main results}

\label{sec:NotionofSolandMainresults}

\subsection{Preliminary definitions and facts}

Throughout the paper, we use classical results of Brenier \cite{Br_87} and Benamou and Brenier \cite{BB_00}, namely, that a transport map exists whenever the source measure is absolutely continuous with respect to Lebesgue-measure, and that it can be characterized in terms of the continuity equation.

\begin{prop}\label{P:transportplan}
	Assume $\mu \in \mcl P_{2,\ac}(\R^d)$ and $\nu \in \mcl P_2(\R^d)$.
	 Then the following hold.
	\begin{enumerate}[(a)]
	\item\label{plan} The infimum in the definition of $d_2(\mu,\nu)$ is attained for a unique measure $\gamma = (\Id, T_\mu^\nu)_\# \mu$, where $T_\mu^\nu: \Rd \to \Rd$ is given by $T_\mu^\nu = \nabla \phi_\mu^\nu$ for a convex function $\phi_\mu^\nu: \RR^d \to \RR$ (which is finite in the support of $\mu$ and unique up to the addition of a constant), and $(T_\mu^\nu)_\#\mu = \nu$.
 	\item\label{BenBre} Whenever $t_0,t_1 \in \R$ with $t_0<t_1$, we have
 \[
		d_2^2(\mu,\nu) =  (t_1 - t_0) \inf \Bigl \{  \int_{t_0}^{t_1} \int_{\RR^d} |\alpha_t(x)|^2 d\mu_t(x) dt \Bigr \},
	\]
 where the infimum is taken over all $\mu \in \mathcal{C}([t_0,t_1], \mathcal{P}_2(\R^d))$ and $\alpha \in L^2([t_0,t_1]\times \R^d, \R^d, dt \otimes \mu_t)$ satisfying, in the sense of distributions,
 \[
		\del_t \mu_t + \div(\alpha_t \mu_t) = 0 \quad \text{in } [t_0,t_1] \times \RR^d, \quad \mu_{t_0} = \mu, \quad \mu_{t_1} = \nu.
	\]
	\end{enumerate}
\end{prop}
We will repeatedly use the notation $T_{\mu}^{\nu}$ for the optimal transport map between the measures $\mu \in \mathcal{P}_{2,ac}(\R^d)$ and $\nu \in \mathcal{P}_2(\R^d)$. 

We will need a regularity result for Brenier's map. Several such results are known, taking advantage of the regularizing properties of the uniformly elliptic Monge-Amp\`ere equation \cite{Caffarelli1991,Cordero2019}. In our setting, we appeal to recent results for optimal transport on unbounded domains, and, in particular, the following special case of \cite[Theorem 1]{Cordero2019}.

\begin{lem}\label{L:bootstrapT}
	Assume $\mu, \nu \in \mcl P_{2,\ac}(\R^d)$ satisfy $\mu, \mu^{-1} ,\nu, \nu^{-1} \in L^\oo_\loc(\R^d)$. Then, for every $R > 0$, there exist $\alpha,\epsilon \in (0,1)$ such that $T_\mu^\nu \in C^\alpha(B_R) \cap W^{1,1+\epsilon}(B_R)$. Moreover, if, for some $k \ge 0$ and $\beta \in (0,1)$, $\mu, \nu \in C^{k,\beta}_\loc$, then $T_\mu^\nu$ is a diffeomorphism of class $C^{k+1,\beta}_\loc$.
\end{lem}

In order to understand the derivatives in the equation \eqref{eq:HJB}, we take advantage of the fact that the $2$-Wasserstein space can be endowed with a Riemannian-like structure. Following \cite{AGS_2008}, for any $\mu \in \mathcal{P}_2(\R^d)$ we define the tangent space of $\mathcal{P}_2(\R^d)$ at $\mu$ as follows.

\begin{defn}\label{D:tangent}
	For any $\mu \in \mcl P_2(\R^d)$,
	\[
		\Tan_\mu(\mcl P_2(\Rd)) := \overline{ \left\{ \nabla \phi : \phi \in C_c^\oo(\Rd) \right\}}^{L^2(\mu;\Rd)}.
	\]
\end{defn}

The following criterion follows from a simple regularization and localization argument; see also \cite[Lemmas D.30,D.31]{FK_book_06} for similar statements.

\begin{lem}\label{L:tangent:crit}
	Assume $\mu \in \mcl P_{2,\ac}(\Rd)$, and $\phi \in L^1_\loc(\R^d)$ is such that $\nabla \phi \in L^1_\loc(\Rd) \cap L^2(\mu)$. Then $\nabla \phi \in \Tan_\mu(\mcl P_2(\Rd))$.
\end{lem}



We now discuss differentiability properties of functionals defined on $\mathcal{P}_2(\R^d)$. It turns out that we will only need such properties at measures which are absolutely continuous with respect to the Lebesgue measure. In the next definition we assume that $F$ is not identically $+\infty$ or $-\infty$.

\begin{defn}\label{D:diff}
	We say that $(r,p) \in \R \times \Tan_{\mu}(\mathcal{P}_2(\R^d))$ belongs to the subdifferential $\partial^{-}F(t,\mu)$ of the functional $F : [0,T] \times \mathcal{P}_2(\R^d) \rightarrow (-\infty, + \infty]$ at $(t,\mu) \in [0,T] \times \mathcal{P}_{2,ac}(\R^d)$ if 
	\begin{align*}
		F(s,\nu) \geq F(t,\mu) +r(s-t) &+ \int_{\R^d} p(x) \cdot (T_{\mu}^{\nu}(x) - x) d\mu(x) + o(d_2(\nu,\mu) + |t-s|)\\
		&\text{as } (s,\nu) \to (t,\mu) \text{ in } [0,T] \times \mathcal{P}_2(\R^d).
	\end{align*}
Similarly, $(r,p) \in \R \times \Tan_{\mu}(\mathcal{P}_2(\R^d))$ belongs to the superdifferential $\partial^{+}F(t,\mu)$ of $F: [0,T] \times \mathcal{P}_2(\R^d) \rightarrow [-\infty, + \infty)$ at $(t,\mu) \in [0,T] \times \mathcal{P}_{2,ac}(\R^d)$ if 
	\begin{align*}
		F(s,\nu) \leq F(t,\mu) +r(s-t) &+ \int_{\R^d} p(x) \cdot (T_{\mu}^{\nu}(x) - x) d\mu(x) + o(d_2(\nu,\mu) + |t-s|)\\
		&\text{as } (s,\nu) \to (t,\mu) \text{ in } [0,T] \times \mathcal{P}_2(\R^d).
	\end{align*}
\end{defn}

\begin{rmk}
	The sub and superdifferential of a functional $F$ on $\mcl P_2(\Rd)$ that is independent of time can be defined analogously by ignoring the time variable. In that case, if $\partial^- F(\mu) \cap \partial^+ F(\mu) = \left \{ p  \right \}$ for some $\mu \in \mcl P_{2,\ac}(\Rd)$, then we say that $F$ is differentiable at $\mu$ and we denote by $\partial_{\mu} F (\mu)= p$ the Wasserstein gradient of $F$.
\end{rmk}

\begin{ex}
    If $f : \R^d \rightarrow \R$ has quadratic growth with continuous derivative of quadratic growth then $\mu \mapsto \int_{\R^d} f(x)d\mu(x)$ is differentiable over $\mathcal{P}_2(\R^d)$ with Wasserstein gradient given, for every $\mu \in \mathcal{P}_2(\R^d)$, by $\nabla f$. If $\nu$ belongs to $\mathcal{P}_2(\R^d)$ then it follows from \cite[Proposition 3.8]{Gigli2004} that $\mu \mapsto d_2^2 (\mu, \nu)$ is differentiable at every $\mu \in \mathcal{P}_{2,ac}(\R^d)$ with Wasserstein gradient given by $2 \bigl(i_d-T_{\mu}^{\nu} \bigr)$.
\end{ex}

We now introduce two well-known functionals that are used in the definition of solution. We define the entropy $\mathcal{E}(\mu)$ of a measure $\mu \in \mathcal{P}_2(\R^d)$ by
\begin{equation}\label{entropy}
	\mcl E(\mu) = 
	\begin{dcases}
		\int_{\Rd} \log \mu(x) d\mu(x) &\text{if } \mu \in \mcl P_{2,\ac}(\Rd),\\
		+\oo & \text{otherwise.}
	\end{dcases}
\end{equation}
As it turns out, even if $\mu (\log \mu)_+ \notin L^1(\Rd)$, the integral defining $\mcl E(\mu)$ is always bounded from below if $\mu \in \mcl P_{2,ac}(\Rd)$. We make this clear with the next result, which also contains a useful coercivity property of the entropy. For $\sigma > 0$, we denote by $g_\sigma \in \mcl P_{\ac}(\Rd)$ the Gaussian measure with variance $\sigma^2$, i.e.
\begin{equation}\label{Gaussian}
	g_\sigma = \frac{1}{(2\pi \sigma)^{d/2} } e^{- \frac{|x|^2}{2\sigma}}dx.
\end{equation}

\begin{lem}\label{L:entpen}
	The functional $\mu \mapsto \mcl E(\mu)$ is lower-semicontinuous with respect to weak convergence of probability measures, and, for any $\sigma > 0$, the inequality
	\begin{equation}\label{entropy:coercive}
		\mcl E(\mu) + \frac{1}{2\sigma} \int_{\Rd} |x|^2 \mu(x)dx \ge - \frac{d}{2} \log(2 \pi \sigma)
	\end{equation}
	holds, with equality if and only if $\mu = g_\sigma$.
\end{lem}

\begin{proof}
If $\mu \in \mcl P_{2,\ac}(\Rd)$ is absolutely continuous with respect to $g_\sigma$ with $\mu = \rho g_\sigma$ for some $\rho \in L^1(g_\sigma)$, then
\begin{equation}\label{entropy:relation}
	\int \rho(x) \log \rho(x) dg_\sigma(x) = \mcl E(\mu) + \frac{1}{2\sigma} \int_{\Rd} |x|^2 \mu(x)dx + \frac{d}{2} \log(2\pi \sigma).
\end{equation}
The left-hand side is the relative entropy of $\mu$ with respect to $g_\sigma$, and thus, by Jensen's inequality, nonnegative, and zero only if $\mu = g_\sigma$. The left-hand side is also lower-semicontinuous with respect to weak convergence (see for instance \cite[Theorem 15.4]{ABS_2021}), and we conclude.
\end{proof}

It is also convenient to introduce the following, nonnegative penalized version of the entropy functional:
\begin{equation}\label{coercive:entropy}
	\mcl E^*(\mu) := \mcl E(\mu) + \pi \int_{\Rd} |x|^2 \mu(x)dx.
\end{equation}
It corresponds to the choice $\sigma = (2\pi)^{-1}$, and, in view of \eqref{entropy:relation}, satisfies
\begin{equation}\label{Estar:coercive}
	\mcl E^* \ge 0 \quad \text{and} \quad \mcl E^*(\mu) \ge \frac{\pi}{2} \int_{\Rd} |x|^2 \mu(x)dx - \frac{d}{2} \log 2.
\end{equation}

The other important functional involved in much of the analysis of the paper is the Fischer information, which is defined as follows:
\begin{equation}
	\mathcal{I}(\mu) :=
	\begin{dcases}
		4\int_{\R^d} \bigl| \nabla \sqrt{ \mu } (x) \bigr|^2 dx &\mbox{ if } \mu \in \mathcal{P}_{2,ac}(\R^d) \mbox{ and } \sqrt{\mu} \in H^1(\R^d), \\
		+ \infty & \mbox{ otherwise.}
	\end{dcases}
\end{equation}
Note that $\mathcal{I}(\mu)$ is finite if and only if $\nabla \log \mu \in L^2(\mu; \R^d)$ (see \cite[Lemma D.7]{FK_book_06}), and then 
$$\mathcal{I}(\mu) = \int_{\R^d} \bigl| \nabla \log \mu(x) \bigr|^2d\mu(x).$$





\begin{lem}\label{lem:FtoE}
	Assume $\mu \in \mcl P_{2,\ac}(\Rd)$ and $\mcl I(\mu) < \oo $. Then $|\mcl E(\mu)| < \oo$.
\end{lem}

\begin{proof}
	The Sobolev embeddings for $\sqrt \mu \in H^1(\Rd)$ in various dimensions imply that
	\[
		\mu \in
		\begin{dcases}
			L^{ \frac{d}{d-2}}(\Rd), & d > 2, \\
			L^p(\Rd), \; 1 \le p < \oo, & d = 2, \\
			L^\oo(\Rd), & d = 1.
		\end{dcases}
	\]
	In any case, we have $\mu (\log \mu)_+ \in L^1(\Rd)$, and so $\mcl E(\mu) < +\oo$. The lower bound $\mcl E(\mu) > -\oo$ follows from \eqref{entropy:coercive} and the fact that $\mu \in \mcl P_2(\Rd)$.
\end{proof}


\subsection{Definition of viscosity solutions and main results}

\label{sec:Defofviscositysolandmainresults}

Before defining our notion of solution for \eqref{eq:HJB}, we make the assumptions on the data precise.

The diffusion matrix satisfies the regularity assumption
\begin{equation}\label{A:a}
	\left\{
	\begin{split}
		&a(\mu,x) = \sigma(\mu,x)\sigma(\mu,x)^*,\\
		&\sigma: \mcl P_2(\R^d) \times \RR^d \to \mcl M_{d \times d} \text{ is uniformly bounded and Lipschitz continuous, and,} \\
		&\text{for all $\mu \in \mcl P_2(\R^d)$, $\sigma(\cdot,\mu) \in C^1(\Rd)$.}	
	\end{split}
	\right.
\end{equation}
We also assume (see Remark \ref{rmk:elliptic} below) that $a$ is uniformly elliptic, that is
\begin{equation}\label{A:aelliptic}
	a(\mu,x)\xi \cdot \xi \ge \kappa|\xi|^2 \text{ for some $\kappa > 0$ and all }\xi \in \RR^d, \; (\mu,x) \in \mcl P_2(\Rd) \times \Rd.
\end{equation}

The Hamiltonian is required to satisfy the following mild regularity and growth assumption:
\begin{equation}\label{AHexistence 3juillet2023}
	\left\{
	\begin{split}
		&\text{for some $C > 0$ and for all $x,y, p,q \in \Rd$ and $\mu,\nu \in \mcl P_2(\Rd)$},\\
		&|H(x,p,\mu) - H(y,q,\nu)| \le C(1 + |p| + |q|)(|x-y| + |p-q| + d_2(\mu,\nu) ).
	\end{split}
	\right.
\end{equation}

We now provide the particular definition of viscosity solution used in this paper. We provide an informal discussion and motivation in the forthcoming subsection.


\begin{defn}\label{D:soln}
   We say that $U : [0,T] \times \mathcal{P}_2(\R^d) \rightarrow \R$ is a viscosity sub-solution of \eqref{eq:HJB} if there is a constant $C>0$ such that, for every $\delta >0$,
	whenever $t \in [0,T]$, $\mathcal{I}(\mu)<\infty$, and $(r,p) \in \del^+ \bigl(U - \delta \mathcal{E} \bigr) (t,\mu)$, we have
\begin{align*}
	-r + \int_{\R^d} \bigl( \div_x a(\mu,x) + a(\mu,x)  \nabla \log \mu(x) \bigr)\cdot p(x) d\mu(x)  + \int_{\Rd} H \bigl(x, p(x), \mu \bigr)d\mu(x) \leq C \delta.
\end{align*}

We say that $V: [0,T] \times \mathcal{P}_2(\R^d) \rightarrow \R$ is a viscosity super-solution if there is a constant $C>0$ such that, for every $\delta >0$, whenever $t \in [0,T]$, $\mathcal{I}(\mu)<\infty$, and $(r,p) \in \del^- \bigl( V + \delta \mathcal{E} \bigr) (t,\mu) $, we have
\begin{align*}
	-r + \int_{\R^d} \bigl( \div_x a(\mu,x) + a(\mu,x)\nabla \log \mu(x) \bigr)\cdot p(x)  d\mu(x) + \int_{\Rd} H \bigl(x, p(x), \mu \bigr)d\mu(x) \geq -C\delta.
\end{align*}
We say that $U$ is a viscosity solution if it is both a sub- and super-solution.
\end{defn}



\begin{rmk}
	By Lemma \ref{lem:FtoE}, the finiteness of the Fischer information of $\mu \in \mcl P_2(\Rd)$ implies that the entropy functional $\mcl E( \mu)$ is well-defined and finite. Furthermore, note that \eqref{AHexistence 3juillet2023} implies that, for some constant $C  > 0$,
	\[
		|H(x,p,\mu)| \le C\left( 1 + |x|^2 + |p|^2 + \int_\Rd |y|^2 d \mu(y) \right) \quad \text{for all } (x,p, \mu) \in \Rd \times \Rd \times \mcl P_2(\Rd).
	\]
	In particular, all terms in the solution inequalities above are well-defined, in view of the fact that $\mu \in \mcl P_2(\Rd)$ and $p \in L^2(\mu)$.
\end{rmk}

\begin{rmk}
	If $\mathcal{I}(\mu)$ is finite and $p \in L^2(\mu) \cap \mathcal{C}^1(\R^d)$, then a formal integration by parts allows the term involving the diffusion matrix $a$ to read
\begin{align*} 
\int_{\R^d} \bigl( \div_x a(\mu,x) + a(\mu,x)\nabla \log \mu(x) \bigr) \cdot  p(x) d\mu(x)= - \int_{\Rd} \tr\left[a(\mu,x) \nabla p(x)\right]d\mu(x),
\end{align*}
as suggested by the equation \eqref{eq:HJB}. However, this computation is only justified if the integral on the right-hand side is itself finite; note that it is not assumed in the definition of solution that $\nabla p \in L^2(\mu)$.
\end{rmk}

\begin{rmk}\label{R:testfunctions}
	An alternative notion of solution involves touching $U$ from above, or $V$ from below, with smooth test functions $\Phi: \mcl P_2(\R^d) \to \RR$ (see \cite{Cosso2021}), which here means that $\Phi$ is differentiable in the sense of Definition \ref{D:diff}, $\del_\mu\Phi$ is continuous on $\mcl P_2(\R^d) \times \Rd$, and $\del_\mu \Phi(\mu, \cdot) \in C^1(\Rd)$ with jointly continuous derivative. This notion is a priori weaker than the one we introduce in this paper, but it is not clear that these two notions of viscosity solution are equivalent, as is the case in the classical finite-dimensional setting of viscosity solutions. We are constrained to work with Definition \ref{D:soln}, because the squared $2$-Wasserstein distance is not a smooth test function in this sense. 
	
\end{rmk}

Our comparison theorem allows to compare sub/super solutions in the class of functions $U:[0,T] \times \mathcal{P}_2(\R^d) \rightarrow \R$ satisfying the following Lipschitz regularity:
\begin{equation}
    U: [0,T] \times \mathcal{P}_2(\R^d) \rightarrow \R \mbox{ is Lipschitz in $\mu \in \mathcal{P}_2(\R^d)$ w.r.t $d_1$ uniformly in $t \in [0,T]$.}
\label{LipregularityAssumption}
\end{equation}

\begin{thm}\label{T:comparison}
	Assume \eqref{A:a} and \eqref{AHexistence 3juillet2023}, and let $U,V: [0,T] \times \mathcal{P}_2(\R^d) \rightarrow \R$ be respectively a bounded sub- and supersolution of \eqref{eq:HJB} satisfying the regularity assumption \eqref{LipregularityAssumption}. Then, for all $(t,\mu) \in [0,T] \times \mcl P_2(\Rd)$,
	\[
		U(t,\mu) - V(t,\mu) \le \sup_{\mcl P_2(\Rd)} \left\{ U(T,\cdot) - V(T,\cdot) \right\}.
	\]
	In particular, given a bounded function $\mcl G$ that is Lipschitz continuous with respect to $d_1$, there is at most one bounded continuous viscosity solution to the terminal value problem \eqref{eq:HJB} satisfying the regularity condition \eqref{LipregularityAssumption}.
\label{UniquenessTheorem3juillet2023}
\end{thm}
The proof of Theorem \ref{UniquenessTheorem3juillet2023} is given in Section \ref{sec:UniquenessSection}.

\begin{rmk}
	The comparison principle still holds when $U$ and $V$ are not necessarily bounded since the $d_1$-Lipschitz assumption implies that they have at most linear growth in the first moment of the measures. We prove Theorem \ref{T:comparison} under the additional boundedness assumption in order to avoid certain technicalities in the proof. 
\end{rmk}


The second main result of the paper is to show that the value function from a mean field control problem gives a viscosity solution of \eqref{eq:HJB} in the sense of Definition \ref{D:soln}. Although more general settings can be considered, we concentrate, for the sake of a simplified presentation, on the case where $a=I_d$ (that is, $\kappa = 1$ in \eqref{A:aelliptic}) and $H(x,p,\mu) = H_1(x,p) - \mathcal{F}(\mu)$. The functions $H_1, \mathcal{F}$ and $\mathcal{G}$ are assumed to satisfy
\begin{equation}\label{AHallinone}
\left \{
\begin{split}
	&\mbox{$H_1 \in \mathcal{C}^2 (\R^d \times \R^d)$, $H_1$ and its derivatives are bounded on $\R^d \times B(0,R)$ for all $R >0$,}\\
	&\mbox{and, for some $C > 0$ and all $(x,p) \in \Rd \times \Rd$,} \\
	&|D_xH_1(x,p)| \leq C (1 +|p|) \mbox{ and}\\
	&\frac{1}{C}I_d \leq D^2_{pp}H_1(x,p) \leq C I_d.
\end{split}
\right.
\end{equation}
Notice that these assumptions imply in particular that $H$ is convex with quadratic growth in the $p$-variable. For example, \eqref{AHallinone} is satisfied when $H_1(x,p) = Ap\cdot p + B(x) \cdot p + h(x)$ for some positive definite symmetric matrix $A \in \mathcal{M}_{d \times d}(\R)$, some $\mathcal{C}^2$ vector field $B : \R^d \mapsto \R^d$ with bounded derivatives and some $\mathcal{C}^2$ function $h : \R^d \to \R$ with bounded derivatives.

For the mean-field costs we assume
\begin{equation}
\mathcal{F}, \mathcal{G}: \mathcal{P}_1(\R^d) \rightarrow \R \mbox{ are $d_1$-Lipschitz continuous and bounded.}
\label{AssumptionsMeanFieldCosts}
\end{equation}
For all $(t_0,m_0) \in [0,T] \times \mathcal{P}_2(\R^d)$, we let
\begin{equation}  \label{def:Uexistence:Sec2}
	U(t_0,m_0) = \inf_{(\mu,\alpha)} \int_{t_0}^T \int_{\R^d} L \bigl(x,\alpha_t(x) \bigr)d\mu_t(x)dt +\int_{t_0}^T \mathcal{F}(\mu_t)dt + \mathcal{G}(\mu_T)
\end{equation}
where $L(x,q) := \sup_{p \in \R^d} \bigl \{ -p\cdot q - H_1(x,p) \bigr \}$ and the infimum is taken over the couples $(\mu,\alpha)$ with $\mu \in \mathcal{C}([t_0,T] ,\mathcal{P}_2(\R^d)), \alpha\in L^2 ([t_0,T] \times \R^d, dt \otimes \mu_t; \R^d) $ satisfying, in the sense of distributions, the Fokker-Planck equation
\begin{equation}
\left \{
\begin{array}{ll}
\partial_t \mu + \div((\alpha_t(x)\mu)-\Delta \mu = 0 &\mbox{ in } (t_0,T) \times \R^d,  \\
\mu(t_0)=\mu_0.
\end{array}
\right.
\label{FPESection2}
\end{equation}
Under these conditions we prove the main existence result.
\begin{thm}
	Assume \eqref{AHallinone} and \eqref{AssumptionsMeanFieldCosts}. Then the value function $U$ defined in \eqref{def:Uexistence:Sec2} is bounded, jointly continuous, and $d_1$-Lipschitz in the measure variable uniformly in time. Moreover, $U$ is a viscosity solution to the HJB equation \eqref{eq:HJB}.
\label{thm:ExistenceThm}
\end{thm}
The proof of Theorem \ref{thm:ExistenceThm} is given in Section \ref{sec:ExistenceSection}. It relies on dynamic programming and requires a careful inspection of the regularity of the maps $t \mapsto d_2(\mu_t,\mu_0)$, $t \mapsto \mathcal{E}(\mu_t)$ and $t \mapsto \mathcal{I}(\mu_t)$ when $\mu_t$ is solution to \eqref{FPESection2} with $\mathcal{I}(\mu_0)<+\infty$ and a regular enough control $\alpha.$

Combining Theorem \eqref{UniquenessTheorem3juillet2023} and Theorem \eqref{thm:ExistenceThm} we get a full well-posedness result.

\begin{cor}
    Assume \eqref{AHallinone} and \eqref{AssumptionsMeanFieldCosts}. Then the value function $U$ is the unique viscosity solution to the HJB equation \eqref{eq:HJB} in the class of bounded continuous functions satisfying the regularity condition \eqref{LipregularityAssumption}.
\end{cor}

\begin{rmk}
  If $\bigl( \Omega, \mathcal{F}, \mathbb{P} \bigr)$ is a probability space supporting a standard $d$-dimensional Brownian motion $(B_t)_{t \geq t_0}$ then, by standard arguments, $U(t_0, \mu_0)$ can be equivalently defined as 
  \begin{equation}
  \label{equivalentformulationU}
  U(t_0, \mu_0) = \inf_{(\alpha_t)_{t \in [t_0,T]}} \int_{t_0}^T \E \left[ L \bigl( X_t, \alpha_t \bigr) \right]dt + \int_{t_0}^T \mathcal{F} \bigl( \mathcal{L}(X_t) \bigr) dt + \mathcal{G}( \bigl( \mathcal{L}(X_T) \bigr), 
  \end{equation}
where $X_t$ solves the stochastic differential equation 
$$X_t = X_0 + \int_{t_0}^t \alpha_t dt + \sqrt{2} \bigl(B_t - B_{t_0} \bigr), \quad t \geq t_0,$$
with $\mathcal{L}(X_0)$, the law of $X_0$, equal to $\mu_0$ and the infimum is taken over the processes $(\alpha_t)_{t \in [t_0,T]}$ in $L^2 \bigl( [t_0,T] \times \Omega; \R^d \bigr)$ adapted to the filtration generated by the Brownian motion. Upon appealing to \cite[Lemma 4.1]{DDJ_23} to approximate $\mathcal{F}$ and $\mathcal{G}$ by differentiable functions, equality \eqref{equivalentformulationU} follows, from the optimality conditions for \eqref{def:Uexistence:Sec2}, see \cite{Briani2018,Daudin2023} for instance, which provide a Lipschitz continuous optimal control for \eqref{def:Uexistence:Sec2} that can be used in the probabilistic formulation and shown, by verification, to be optimal. Notice that equality \eqref{equivalentformulationU} also follows from our comparison principle since both formulations yield a  viscosity solution to the HJB equation.
\label{EquivalenceProbabilityformulation}
\end{rmk}

\subsection{Motivation for the definition of viscosity solution}
Consider the model case
\begin{equation}\label{eq:HJBidea}
    \left \{
	\begin{split}
	&-\partial_t U(t,\mu) - \kappa \int_{\R^d} \div_x\partial_{\mu}U(t,\mu,x)d\mu(x)\\
	&\qquad +  \int_{\R^d} h\left( \partial_{\mu}U(t,\mu,x) \right) d\mu(x) = 0 \quad \mbox{in }[0,T] \times \mathcal{P}_2(\R^d) \\
	&U(T,\mu) = \mathcal{G}(\mu) \quad \mbox{in } \mathcal{P}_2(\R^d),
	\end{split}
    \right.
\end{equation}
that is, $a= \kappa I_d$ for some $\kappa > 0$ and $H(x,p,\mu) = h(p)$ for some locally Lipschitz function $h : \R^d \rightarrow \R$.

Our approach to proving the comparison principle is to double variables, as in the classical setting of viscosity solutions, using the squared $2$-Wasserstein distance as a penalizing function. Unlike in the standard finite-dimensional setting, this quadratic penalizing function is not smooth on all of $\mcl P_2(\Rd) \times \mcl P_2(\Rd)$, and we must introduce further penalizations in order to force the extrema to be attained at sufficiently ``nice'' measures. 

In particular, for a sub and supersolution $U$ and $V$ and small parameters $\delta,\epsilon,\gamma > 0$, we consider the maximum point\footnote{We actually use the modified functional $\mcl E^*$ defined above in \eqref{coercive:entropy} for added coercivity, but we ignore such technicalities here in order to keep the discussion simple.} $(\bar s, \bar t, \bar \mu, \bar \nu)$ of 
\[
	(s,t,\mu,\nu) \mapsto U(s,\mu) - V(t,\nu) - \frac{1}{2\epsilon}d_2^2(\mu,\nu) - \delta \mcl E(\mu) - \delta \mcl E(\nu) - \frac{1}{2\epsilon}|s-t|^2 - \gamma t.
\]
One of the key results of the paper, Proposition \ref{P:JKO}, states that, if $U$ and $V$ are Lipschitz continuous with respect to $d_1$, then any such maximum must occur for absolutely continuous measures $\bar \mu$ and $\bar \nu$ that are locally Lipschitz-continuous, have full support, and whose Fischer informations $\mcl I(\bar \mu)$ and $\mcl I(\bar \nu)$ are finite. In particular, the regularity theory for the Monge-Amp\`ere equation implies that $T_{\bar{\mu}}^{\bar{\nu}}$ is a $\mathcal{C}^1$ diffeomorphism. Differentiability properties of the $2$-Wasserstein distance, see \cite[Proposition 3.8]{Gigli2004}, lead to 
\[
	\del_\mu d_2^2(\cdot,\bar \nu) [\bar \mu] = i_d - T_{\bar \mu}^{\bar \nu}
	\quad \text{and} \quad
	\del_\nu d_2^2(\bar \mu, \cdot) [\bar \nu] = i_d - T_{\bar \nu}^{\bar \mu}.
\]
while a formal computation implies
\begin{equation}\label{entropy:derivative}
	\del_\mu \mcl E(\mu,x) = \nabla \log \mu(x).
\end{equation}
If $U$ and $V$ were smooth, this would lead to the solution inequalities
\begin{equation}\label{eq:subsolforUintro2023}
	\begin{split}
	&\frac{\bar t - \bar s}{\epsilon}  - \frac{\kappa}{\epsilon} \int_{\R^d} \tr \Bigl( I_d - \nabla T_{\bar{\mu}}^{\bar{\nu}}(x) \Bigr) d\bar{\mu}(x) - \delta \kappa \int_\Rd \Delta \log \bar \mu(x) d \bar \mu(x)\\
	&+ \int_{\R^d} h\left( \frac{ x - T_{\bar{\mu}}^{\bar{\nu}}(x)}{\epsilon} + \delta \nabla \log \bar \mu(x) \right)  d\bar{\mu}(x) \le 0
	\end{split}
\end{equation}
and
\begin{equation}\label{eq:supersolforVintro2023}
	\begin{split}
	&\frac{\bar t - \bar s}{\epsilon} - \gamma  +  \frac{\kappa}{\epsilon}\int_{\R^d} \tr \Bigl( I_d - \nabla T_{\bar{\nu}}^{\bar{\mu}}(y) \Bigr) d\bar{\nu}(y)  + \delta \kappa \int_\Rd \Delta \log \bar \nu(y) d\bar \nu(y)\\
	&+  \int_{\R^d} h\left( \frac{T_{\bar{\nu}}^{\bar{\mu}}(y) - y}{\epsilon} - \delta \nabla \log \bar \nu(y) \right) d\bar{\nu}(y) \geq 0.
	\end{split}
\end{equation}
Since $U$ and $V$ are $d_1$-Lipschitz and $\bar{\mu}$, $\bar{\nu}$ have full support, the technical result Proposition \ref{P:JKO} implies that
$$ \norm{\frac{1}{\epsilon} \bigl(i_d - T_{\bar{\mu}}^{\bar{\nu}} \bigr) + \delta \nabla \log \bar{\mu}}_{\infty} \leq L, \quad  \norm{\frac{1}{\epsilon} \bigl(i_d - T_{\bar{\nu}}^{\bar{\mu}} \bigr) - \delta \nabla \log \bar{\nu}}_{\infty} \leq L  $$
where $L$ is a common Lipschitz constant for $U$ and $V$. Denote by $\bar \pi = (i_d, T_{\bar \mu}^{\bar \nu})_\# \bar \mu$ the optimal transport plan between $\bar \mu$ and $\bar \nu$. Subtracting \eqref{eq:supersolforVintro2023} from \eqref{eq:subsolforUintro2023} thus leads, in view of the local Lipschitz continuity of $h$, for some $C_{h,L}>0$ depending only on the Lipschitz constant of $h$ in a ball of radius $L$ centered at the origin,
\begin{equation}\label{laststep}
\begin{split}
	\gamma &- \delta \kappa \int_\Rd \Delta \log \bar{\mu}(x) \bar \mu(x)dx
	- \delta \kappa \int_\Rd \Delta \log \bar{\nu}(x) \bar \nu(x)dx\\
	&\le \int_{\Rd\times \Rd} \left[ h\left( \frac{x-y}{\epsilon} - \delta \nabla \log \bar \nu(y)\right) - h\left( \frac{x-y}{\epsilon} + \delta \nabla \log \bar \mu(x) \right) \right] d\bar \pi(x,y) \\
	&+ \frac{\kappa}{\epsilon} \int_{\Rd } \tr\left[ 2I_d - \nabla T_{\bar \mu}^{\bar \nu}(x) - \nabla T_{\bar \mu}^{\bar \nu}(x)^{-1} \right] d\bar{\mu}(x) \\
	&\le C_{h,L} \delta\int_\Rd \left( \left| \nabla \log \bar{\mu}(x) \right| \bar \mu(x) + \left| \nabla \log \bar{\nu}(x) \right| \bar \nu(x)  \right)dx,
\end{split}
\end{equation}
where the last line uses the simple inequality $A + A^{-1} \geq 2I_d$ for any positive definite matrix $A$; recall that $\nabla T_\mu^\nu = \nabla^2 \phi_\mu^\nu$ for the convex function $\phi_\mu^\nu$ from Proposition \ref{P:transportplan}.

Integrating by parts on the left-most side of \eqref{laststep}, we find that
\begin{equation}\label{absorb}
	\gamma \le - \delta \int_{\Rd}\left( \kappa |\nabla \log \bar \mu|^2  - C_{h,L} |\nabla \log \bar \mu| \right)d \bar \mu - \delta \int_{\Rd}\left( \kappa |\nabla \log \bar \nu|^2  - C_{h,L}|\nabla \log \bar \nu| \right)d \bar \nu,
\end{equation}
and we conclude that
\begin{equation}\label{noticeCdependsonkappa}
	\gamma \le \frac{ C_{h,L}}{2\kappa} \delta,
\end{equation}
which is a contradiction for sufficiently small $\delta$.

As it turns out, the formal identity \eqref{entropy:derivative} cannot be easily verified, at least not in the sense of Definition \ref{D:diff} above. Therefore, the definition of solution above is formulated in terms of elements of the sub and superdifferential of respectively $U - \delta \mcl E$ and $V + \delta \mcl E$, with a small error term proportional to $\delta$ appearing on the right-hand side of the solution inequalities.

\begin{rmk}[Uniform ellipticity] \label{rmk:elliptic}
It is crucial in the formal argument above, especially in the final step \eqref{laststep}, that the second-order term in the equation \eqref{eq:HJB} be \emph{semilinear}, which allows for the integration by parts. 
Notice also that the constant appearing in front of $\delta$ on the right-hand side of \eqref{noticeCdependsonkappa} depends inversely on the ellipticity $\kappa$. Therefore, in order for Definition \ref{D:soln} to be properly motivated, we will always take the uniform ellipticity assumption in what follows. This absorbing of the terms must then also play a role in any proof of existence, as in Section \ref{sec:ExistenceSection} for the value function from mean field control.

Let us comment further on the method we use to bound \eqref{absorb}. On the one hand, even if $\kappa I_d$ is replaced with a degenerate elliptic matrix $a$ (or even if $\kappa = 0$), the terms involving $|\nabla \log \bar\mu|^2$ and $|\nabla \log \bar\nu|^2$ have a beneficial sign. However, although the terms $\delta \int_\Rd |\nabla \log \bar \mu| d\bar \mu$ and $\delta \int_\Rd |\nabla \log \bar \nu |d \bar \nu$ turn out to be bounded uniformly in $\delta$, they cannot be expected in general to vanish with $\delta$. We illustrate this with a simple example: set $\phi(x) := \min (1, |x|)$ and define $U(\mu) := -\int \phi d \mu$, which is $d_1$-Lipschitz with constant $1$. For $\delta > 0$, the function
	\[
		U(\mu) - \delta \mcl E(\mu) - \delta \int |x|^2 \mu(dx) = \int_\Rd \left( -\phi(x) - \delta \log \mu(x)  - \delta |x|^2 \right)\mu(x) dx
	\]
	is maximized over $\mcl P_2(\Rd)$ by $\bar \mu(x) := \bar c \exp\left( - \frac{1}{\delta} \phi(x) - |x|^2 \right)$, where $\bar c$ is a normalizing constant. We then compute
	\[
		\delta \int_\Rd |\nabla \log \bar \mu(x)| d \bar\mu(x) = \int_\Rd | \nabla \phi(x) + 2\delta x| d \bar \mu(x).
	\]
	As $\delta \to 0$, $\bar \mu$ converges to $\delta_0$, which is the unique maximizer for $U$. In a neighborhood of $0$, $|\nabla \phi(x) + 2 \delta x| = (1 + 4 \delta |x| + 4\delta^2|x|^2)^{1/2}$, and so
	\[
		\lim_{\delta \to 0} \delta \int_\Rd |\nabla \log \bar \mu(x)| d \mu(x) = 1 \ne 0.
	\]
\end{rmk}

\begin{rmk}[Hilbertian Approach]
\label{rmk:HilbertianApproach}

We comment briefly on the approach put forward by P.-L. Lions \cite{PLLCdF} in lectures given at the Coll\`ege de France for converting an equation set on $\mcl P_2(\R^d)$ to one set on a Hilbert space. This involves ``lifting'' a regular function $U: \mcl P_2(\R^d) \to \RR$ to a function $u$ defined on the Hilbert space
$\mcl H := L^2(\Omega; \RR^d)$ by setting
\begin{equation}\label{lift}
	u(X) := U(\mcl L_X), \quad X \in \mcl H,
\end{equation}
where $\mcl L_X \in \mcl P_2(\RR^d)$ denotes the law of $X$. 
The merit of this approach is that the nonlinear space $\mcl P(\R^d)$ is replaced by the separable Hilbert space $\mcl H$, which raises the potential for treating the relevant equations with the theory of viscosity solutions on infinite dimensional spaces \cite{Crandall1992,PLL1988,PLL1989,CL_85,CL_86}. In particular, one has access to a large family of smooth test functions that do
not suffer from the same subtle regularity issues as the squared $2$-Wasserstein distance on $\mcl P_2(\R^d)$. This method has been applied to the study of many mean field control and mean field game equations in settings without noise, leading to first-order equations,
or with common noise, in which case the second-order Fr\'echet derivatives over $\mcl H$ arise; see for instance \cite{CG_19,CS_22,GT_19, MS_23}.

In the setting of the current paper, the second-order term in \eqref{eq:HJB} is expected to describe mean field problems affected by idiosyncratic, rather than common, noise, and so the lift of this term to the Hilbert space can be formally understood by taking advantage of independence in the probability space $\Omega$. If, for some $X \in \mcl H$, there exists a random variable $Z \in \mcl H$ such that
\begin{equation}\label{Z?}
	\E[Z] = 0, \quad \E[Z \otimes Z] = I_d, \quad \text{and} \quad \text{$Z$ is independent of $X$},
\end{equation}
then 
\[
	\nabla^2_Xu(X)[Z,Z] = \int_{\Rd}  \div_x \del_\mu U(\mcl L_X, x)\mcl L_X(dx),
\]
which matches the second-order term in the model equation \eqref{eq:HJBidea}.
However, in applying this idea to viscosity solution arguments,
the main obstruction, one which seems to have been overlooked in previous contributions to this subject, is that, for a given $X \in \mcl H$, a random variable $Z \in \mcl H$ satisfying \eqref{Z?} need not exist in general. Indeed, consider the probability space $\Omega = [0,1]^d$ with the Lebesgue measure, and let $X \in L^2([0,1]^d; \RR^d)$ be the function $X(x) = x$. Then clearly \emph{no} nontrivial random variable $Z \in L^2([0,1]^d; \Rd)$ is independent of $X$, since we can write $Z(x) = Z(X(x))$ for $x \in [0,1]^d$. Moreover, expanding the probability space to allow for more independence completely changes the nature of the space $\mcl H$, which interferes with perturbative arguments that are commonly used in the theory of viscosity solutions.
\end{rmk}

\subsection{Examples}

\label{sec:exmaples}

We conclude this section with two instructive examples in which the equation \eqref{eq:HJB} reduces to equations on finite-dimensional spaces. We consider in particular the model case \eqref{eq:HJBidea} where $h(p) := \frac{1}{2}|p|^2$ and $a= I_d$. 

When $\mathcal{G}(\mu) = \int_{\R^d} g_1(x)d\mu(x)$ for some bounded, Lipschitz continuous map $g_1: \R^d \rightarrow \R$, the unique viscosity solution to \eqref{eq:modelcase} is given by $U(t,\mu) = \int_{\R^d} u_1(t,x)d\mu(x)$, where $u_1$ is (the classical) solution to 
\begin{equation}\label{eq:viscousHJBfinitedim}
    \left \{
    \begin{array}{ll}
\displaystyle    -\partial_t u_1 - \Delta u_1 +  \frac{1}{2} |\nabla u_1|^2 = 0 & \mbox{in }[0,T] \times \R^d, \\
\displaystyle    u_1(T,x) = g_1(x) & \mbox{in } \R^d.
    \end{array}
    \right.
\end{equation}
This example indicates in particular that, in view of our proof of Theorem \ref{T:comparison}, the comparison principle for the second-order equation \eqref{eq:viscousHJBfinitedim} can be proved without appealing to the lemma of Ishii and the technology developed for second order, finite-dimensional, equations, see \cite{Crandall1992}. Of course, this observation is dependent on the fact that the solution of \eqref{eq:viscousHJBfinitedim} is smooth, as a consequence of the nondegenerate ellipticity of the equation. 

When $\mathcal{G}(\mu) = g_2 \bigl( \langle \mu \rangle \bigr)$ for some bounded, Lipschitz continuous map $g_2:\R^d \rightarrow \R$, where $\langle \mu \rangle = \int_{\R^d} x d\mu(x)$ is the expectation of $\mu$, then the unique viscosity solution 
to \eqref{eq:modelcase} is given by $U(t,\mu) = u_2( t, \langle \mu \rangle )$ where $u_2$ is the unique viscosity solution (in the classical sense) to 
\begin{equation}\label{eq:modelcase}
    \left \{
    \begin{array}{ll}
\displaystyle    -\partial_t u_2 +  \frac{1}{2} |\nabla u_2|^2 = 0 & \mbox{in }[0,T] \times \R^d, \\
\displaystyle    u_2(T,x) = g_2(x) & \mbox{in } \R^d.
    \end{array}
    \right.
\end{equation}
This is a consequence of our well-posedness result, which gives a control formulation for the solution, as well as Proposition 2.10 in \cite{DDJ_23}.  This second example shows, in particular, that no regularizing effect with respect to the measure argument can be expected in general, even if the diffusive term is non-degenerate.

\section{Entropy regularization}

\label{sec:EntropyRegularization}



The main result of this section is the following proposition. Recall that the notation $\mathcal{E}^*$ was introduced in \eqref{coercive:entropy}.

\begin{prop}\label{P:JKO}
Assume that $U,V: \mcl P_2(\R^d) \to \R$ are Lipschitz continuous with respect to $d_1$, and fix $\epsilon,\delta > 0$. Then
\begin{align}
    \mathcal{P}_2(\R^d) \times \mathcal{P}_2(\R^d) \ni (\mu,\nu) &\mapsto U(\mu) - V(\nu) - \frac{1}{2\epsilon} d_2^2(\mu,\nu) - \delta \mathcal{E}^*(\mu) - \delta \mathcal{E}^*(\nu) 
\label{MaximumDoublingofVariables}
\end{align}
achieves its supremum. Any maximum point $(\bar{\mu}, \bar{\nu})$ belongs to $ \mathcal{P}_{2,\ac}(\R^d) \times \mathcal{P}_{2,\ac}(\R^d)$ and satisfies the following properties:
\begin{enumerate}[(a)]
    \item\label{measure:bounds} $\bar{\mu}, \bar{\nu}, \bar{\mu}^{-1}, \bar{\nu}^{-1} \in W^{1,\oo}_{\loc}(\R^d)$, $\mathcal{I}(\oline \mu)<+\infty$, $\mathcal{I}(\oline \nu)<+\infty$ and $T_{\bar{\mu}}^{\bar{\nu}}, T_{\bar{\nu}}^{\bar{\mu}} \in \mathcal{C}^{1,\beta}_{loc}$ for all $\beta \in (0,1)$.

    \item\label{diff:bounds} Define
        \[
            p(x) := \frac{1}{\epsilon} \bigl(x-T_{\bar{\mu}}^{\bar{\nu}}(x) \bigr) +2\pi \delta x, \quad x \in \Rd
        \]
        and
        \[
        		q(y) := \frac{1}{\epsilon}\bigl(T_{\bar{\nu}}^{\bar{\mu}}(y)-y \bigr) - 2\pi \delta y, \quad y \in \Rd.
        \]
        Then $p \in \del^+ \bigl(U - \delta \mathcal{E} \bigr) (\oline \mu)$, $q \in \del^- \bigl(V + \delta \mathcal{E} \bigr) (\oline \nu)$ and
        \begin{equation}\label{boundedp}
        		\norm{p + \delta \nabla \log \bar{\mu}}_{L^\oo} \leq \Lip(U;d_1), \quad \norm{q - \delta \nabla \log \bar{\nu}}_{L^\oo} \leq \Lip(V, d_1).
	\end{equation}

\end{enumerate}
\end{prop}


The proof of Proposition \ref{MaximumDoublingofVariables} requires several preliminary facts.

For $F: \mcl P_2(\R^d) \to \R$, we introduce a semicontinuity property along lines in $\mcl P_2(\R^d)$:
\begin{equation}\label{A:semicty}
	\left\{
	\begin{split}
	&\text{for all $R > 0$, there exists $M_R > 0$ such that,}\\
	&\text{if $\mu,\tilde \mu \in \mcl P_2(\R^d)$ and }\mcl M_2(\mu) + \mcl M_2(\tilde \mu) \le R, \text{ then}\\
	&\limsup_{\lambda \to 0^+} \frac{ F(\lambda \tilde \mu + (1-\lambda)\mu) - F(\mu) }{\lambda} \le M_R .
	\end{split}
	\right.
\end{equation}
Recall that $\mathcal{M}_2(\nu) = \int_{\R^d} |x|^2d\nu(x)$ denotes the second order moment of the measure $\nu \in \mathcal{P}_2(\R^d)$. 

\begin{lem}\label{lem:semicty}
The condition \eqref{A:semicty} is satisfied by $F$ if
\begin{enumerate}[(a)]
	\item\label{Lipsemicty} $F$ is $d_1$-Lipschitz continuous, or
	
	\item\label{d22semicty} $F = d_2^2(\cdot,\nu)$ for any fixed $\nu \in \mcl P_2(\Rd)$.
\end{enumerate}
\end{lem}

\begin{proof}
	Let $\mu$, $\tilde \mu$ be as in \eqref{A:semicty}. If $F$ is $d_1$-Lipschitz, then part \eqref{Lipsemicty} is proved upon estimating
\[
	F(\lambda \tilde \mu + (1-\lambda) \mu) - F(\mu) \le \lambda \sup_{\norm{\nabla f}_\oo \le 1} \int_\Rd f(x) d(\tilde{\mu} - \mu)(x)
	\le \lambda\left( \mcl M_2(\tilde{\mu})^{1/2} + \mcl M_2(\mu)^{1/2} \right).
\]
We establish part \eqref{d22semicty} upon computing, by convexity,
\[
	d^2_2(\lambda \tilde \mu + (1-\lambda) \mu,\nu) - d^2_2(\mu, \nu) \le \lambda\left( d^2_2(\tilde \mu,\nu) - d^2_2(\mu,\nu) \right).
\]
\end{proof}

\begin{lem}\label{L:entropyreg}
	Assume that $F: \mcl P(\R^d) \to \R$ satisfies \eqref{A:semicty}, and suppose that
	\[
		\mcl P_2(\R^d) \ni \mu \mapsto F(\mu) + \mcl E(\mu)
	\]
	is bounded from below and attains a minimum at $\bar \mu \in \mcl P_2(\R^d)$. Then $\bar \mu \in \mcl P_{2,\ac}(\R^d)$ and $\bar \mu^{-1} \in L^\oo_\loc$.
\end{lem}

\begin{proof}
	Fix $r > 0$ and $x_0 \in B_r \subset \R^d$, and, for $\delta \in (0,1)$, set $\nu_\delta(x) := g_\delta(x - x_0)$, where $g_\delta$ is the Gaussian measure \eqref{Gaussian}. Note that $M_2(\nu_\delta)$ is bounded independently of $\delta$, while \eqref{entropy:relation} implies that $|\mcl E(\nu_\delta)| < +\oo$ for all $\delta \in (0,1)$.
		
	 For $\lambda \in (0,1)$, define $\nu_{\delta,\lambda} = (1-\lambda)\bar{\mu} + \lambda \nu_{\delta}$. Then $F(\nu_{\delta,\lambda}) + \mcl E(\nu_{\delta,\lambda}) \ge F(\bar \mu) + \mcl E(\bar \mu)$, and so \eqref{A:semicty} gives, for some $R > 0$ depending only on $r$ and $M_2(\bar \mu)$, for all $\delta,\lambda \in (0,1)$,
	 \[
	 	\mcl E(\bar \mu) - \mcl E(\nu_{\delta,\lambda})  \le M_R \lambda.
	 \]
	 By convexity of $s \mapsto s \log s$,
	 \begin{align*}
	 	\mcl E(\bar \mu) - \mcl E(\nu_{\delta,\lambda})
		\ge \int_{\Rd} (1 + \log \nu_{\delta,\lambda}(x)) (\bar \mu(x) - \nu_{\delta,\lambda}(x))dx
		= \lambda \int_{\Rd} (1 + \log \nu_{\delta,\lambda}(x)) (\bar \mu(x) - \nu_\delta(x))dx.
	 \end{align*}
	 We then find that
	 \[
	 	\int_{\Rd}  \log \nu_{\delta,\lambda}(x) (\bar \mu(x) - \nu_\delta(x))dx
		\le M_R.
	 \]
	 Since $\log$ is increasing, we have (distinguishing between the points $x$ where $\bar{\mu}(x) > \nu_{\delta}(x)$ and those where $\bar{\mu}(x) \leq \nu_{\delta}(x)$)
	 \[
	 	\log \nu_{\delta,\lambda}(x) (\bar \mu(x) - \nu_\delta(x)) \ge \log \nu_\delta(x) (\bar \mu(x) - \nu_\delta(x)) = -\left( \frac{|x|^2}{2\delta} + \frac{d}{2} \log(2\pi \delta) \right) (\nu_\delta(x) - \bar \mu(x)).
	 \]
	 The right-hand side belongs to $L^1(\Rd)$ because $\bar \mu, \nu_\delta \in \mcl P_{2,\ac}(\R^d)$, and therefore we may send $\lambda \to 0^+$ and invoke Fatou's lemma to infer that
	 \[
	 	\int_{\Rd} \log \bar \mu d(\bar \mu - \nu_\delta) \le M_R.
	 \]
	 We thus have
	\[
	 	- (\log \bar \mu * g_\delta)(x_0) \le M_R - \mcl E(\bar \mu),
	 \]
	 and so conclude that $\bar \mu$ is uniformly bounded below on $B_r$ upon sending $\delta \to 0$.
	 
\end{proof}


\begin{lem}\label{lem:boundedp}
	Assume that $\mathrm{Lip}(U;d_1) \le L$, $\mu_0 \in \mcl P_{2,\ac}(\Rd)$ is such that $\mcl E(\mu_0) < \oo$, $\delta \ge 0$, and $p \in \del^+(U - \delta \mcl E)(\mu_0)$ (resp. $p \in \del^-(U + \delta \mcl E)(\mu_0)$). Then $|p +  \delta \nabla \log \mu_0| \le L$ (resp. $|p - \delta \nabla \log \mu_0| \le L$) a.e. in the support of $\mu_0$.
\end{lem}

\begin{rmk}
	When $\delta = 0$, the requirement that $\mcl E(\mu_0) < +\oo$ can be dropped.
\end{rmk}

\begin{proof}[Proof of Lemma \ref{lem:boundedp}]
	We prove only the superdifferential statement, since the proof is analogous for the subdifferential. Let $v: \RR^d \to \RR^d$ be smooth and bounded, and consider, for some $t_0 \geq 0$, the flow
	\begin{equation}\label{v:cty}
		\del_t \mu_t + \div(v(x) \mu_t) = 0 \quad \text{in } (t_0,\oo), \quad \mu|_{t = t_0} = \mu_0.
	\end{equation}
	In view of the smoothness of $v$, the map $[t_0,T] \ni t \mapsto \mu_t \in L^1(\Rd)$ is continuous. 
	
	The Benamou-Brenier formula (Proposition \ref{P:transportplan}\eqref{BenBre}) allows for the bound
	\begin{equation}\label{d2linear}
		d_2(\mu_t , \mu_0) \le \norm{v}_\oo |t-t_0|.
	\end{equation}
	We also estimate
	\begin{align*}
		d_1(\mu_t,\mu_0) &= \sup_{\norm{\nabla f}_\oo \le 1} \int_{\RR^d} f(x) d(\mu_t - \mu_0)(x)\\
		&= \sup_{\norm{\nabla f}_\oo \le 1} \int_{t_0}^t \int_{\RR^d} \nabla f(x) \cdot v(x)d\mu_s(x)ds
		\le \int_{t_0}^t \int_{\R^d}|v(x)| d \mu_s(x) ds,
	\end{align*}
	and so
	\[
		|U(\mu_t) - U(\mu_0)| \le L \int_{t_0}^t \int_{\R^d} |v(x)| d \mu_s(x) ds.
	\]
	
	The vector field $v$ being smooth, \eqref{v:cty} satisfies the renormalization property, and so, for any smooth and bounded $e: \R_+ \to \R$, we may compute
	\[
		\del_t e(\mu_t) + \div(v(x) e(\mu_t)) = (\div v(x)) ( e(\mu_t) - \mu_t e'(\mu_t)).
	\]
	By standard regularization and localization procedures as in for instance \cite{diperna_lions_89}, we may approximate the function $s \mapsto s \log s$ by such functions $e$, using the fact that $\mcl E(\mu_0) < +\oo$ and $\mu_t \in \mcl P_2(\Rd)$ for $t > t_0$, to deduce that $\mcl E(\mu_t) < +\oo$ for $t > t_0$ and 
	\[
		\mcl E(\mu_t) - \mcl E(\mu_0) = -\int_{t_0}^t \int_\Rd \div v(x) d\mu_s(x)ds.
	\]
Fix $\varphi \in C_c^\oo(\Rd)$ and $h>0$. Introducing the transport plan
	\[
		\eta_h = \Bigl(i_d, T_{\mu_0}^{\mu_{t_0+h}} -i_d \Bigr)_\# \mu_0,
	\]
	we can write
	\begin{align}
	\notag	 \int_{\R^d} \varphi(x) d\mu_{t_0+h}(x) - \int_{\R^d} \varphi(x) d\mu_0(x)  
	\notag	&=  \int_{\R^d \times \R^d} \bigl[ \varphi(x+y) - \varphi(x) \bigr]d\eta_h(x,y)\\
	\notag	&= \int_{\R^d \times \R^d} \nabla \varphi(x) \cdot y d\eta_h(x,y) + R_h\\
		&= \int_{\R^d} \nabla \varphi(x) \cdot \left( T_{\mu_0}^{\mu_{t_0+h}}(x) - x \right) d\mu_0(x) + R_h,
  \label{tobeusedlaterinSection5}
	\end{align}
	where
	\begin{align}
	\notag	|R_h| &\leq  \norm{\nabla^2 \phi}_{L^\oo} \int_{\Rd \times \Rd} |y|^2 d\eta_h(x,y)\\
			&= \norm{\nabla^2 \phi}_{L^\oo}  \int_{\R^d}  \bigl| T_{\mu_0}^{\mu_{t_0+h}}(x) -x \bigr|^2 d\mu_0(x) =    \norm{\nabla^2 \phi}_{L^\oo} d_2^2( \mu_{t_0+h}, \mu_0).
   \label{tobeusedlateraswell}
	\end{align}
 On the other hand, using the equation satisfied by $\mu$ we have that 
 $$ \int_{\R^d} \varphi(x) d\mu_{t_0+h}(x) - \int_{\R^d} \varphi(x) d\mu_0(x)  = \int_{t_0}^{t_0+h} \int_{\R^d}  \nabla \varphi(x) \cdot v(x)d\mu_s(x)ds$$
 and therefore
	\begin{align*}
		&\left| \int_\Rd \nabla \varphi(x) \cdot (T_{\mu_0}^{\mu_{t_0+h}}(x) - x)d\mu_0(x)
		- \int_{t_0}^{t_0+h} \int_{\Rd}  \nabla \varphi(x) \cdot v(x) d\mu_s(x)ds \right| \\
		& \le \norm{\nabla^2 \varphi}_\oo d^2_2(\mu_{t_0+h},\mu_0) .
	\end{align*}
	Fix $\epsilon > 0$. It then follows from the superdifferential property, the Cauchy-Schwarz inequality, and the estimate \eqref{d2linear} that, for all sufficiently small $h$ depending on $\norm{v}_\oo$ and $\epsilon$,
	\begin{align*}
		&\int_{t_0}^{t_0+h} \int_\Rd \left( \nabla \varphi(x) \cdot v(x) - \delta \div v(x) \right) d \mu_s(x) \\
		&\ge  - L \int_{t_0}^{t_0+h} \int_{\R^d} |v(x)| d \mu_s(x) ds - \norm{\nabla^2 \varphi}_\oo \norm{v}_\oo^2 h^2 - \norm{\nabla \varphi - p}_{L^2(\mu_0)} \norm{v}_\oo h - \epsilon h.
	\end{align*}
	We now divide by $h$ and send $h \to 0^+$. Then, since $v$ and $\varphi$ are smooth and bounded and $\mu_s \to \mu_0$ as $s \to t_0^+$ in $L^1(\Rd)$ this leads to,
	\begin{align*}
		 \int_\Rd \left( \nabla \varphi(x) \cdot v(x)  - \delta \div v(x) \right)d \mu_0(x)
		\ge - L \int_\Rd |v(x)| d \mu_0(x)  - \norm{\nabla \varphi - p}_{L^2(\mu_0)} \norm{v}_\oo - \epsilon.
	\end{align*}
	Since $\epsilon$ and $\varphi$ were arbitrary and $p \in \Tan_{\mu_0}(\mcl P_2(\Rd))$, we infer that
	\[
		\int_\Rd \left( p(x) \cdot v(x) - \delta \div v(x) \right)d \mu_0(x)
		\ge - L \int_\Rd |v(x)| d \mu_0(x).
	\]
	In other words, since $v$ was arbitrary, we may conclude that the distribution
	\[
		p + \delta \nabla \log \mu_0 = p + \delta \frac{\nabla \mu_0}{\mu_0} 
	\]
	is bounded on $L^1(\mu_0)$ with operator norm $L$, and so, in particular, belongs to $L^\oo(\mu_0)$ with $\norm{ p + \delta \nabla \log \mu_0}_{L^\oo(\mu_0)} \le L$.
	\end{proof}


We now present the proof of Proposition \ref{P:JKO}. 
%

\begin{proof}[Proof of Proposition \ref{P:JKO}]
Define
\[
    \Phi(\mu,\nu) = U(\mu) - V(\nu) - \frac{1}{2\epsilon} d_2^2(\mu,\nu) - \delta \mathcal{E}^*(\mu) - \delta \mathcal{E}^*(\nu).
\]
In view of \eqref{Estar:coercive},
\[
	\delta \mathcal{E}^*(\mu)  
	\ge \frac{\delta \pi}{2} \int_{\Rd} |x|^2 d\mu - \frac{\delta d}{2} \log(2),
\]
with equality achieved for the Gaussian measure $\mu = g_{1/\pi}$ (recall the notation from \eqref{Gaussian}), while the $d_1$-Lipschitz continuity of $U$ implies
\[
	U(\mu) \le U(\delta_0) + \Lip(U,d_1)d_1(\mu,\delta_0) = U(\delta_0) + \Lip(U,d_1) \int_\Rd |x| \mu(dx).
\]
Similar bounds hold for the terms involving $\nu$, and we therefore see that $\Phi$ is bounded from above. 

For $n = 1,2,\ldots$, let $(\mu_n,\nu_n) \in \mcl P_{2,\ac}(\R^d) \times \mcl P_{2,\ac}(\R^d)$ be such that $\Phi(\mu_n,\nu_n) \ge \sup \Phi - \frac{1}{n}$. From the inequality $\Phi(g_{1/\pi}, g_{1/\pi}) \le \Phi(\mu_n,\nu_n) + \frac{1}{n}$, we deduce that
\begin{align*}
	\frac{1}{2\epsilon} d^2_2(\mu_n,\nu_n) &+ \frac{\pi \delta}{2} \int_{\Rd} |x|^2 d\mu_n(x) + \frac{\pi \delta }{2}\int_{\Rd} |x|^2 d\nu_n(x) \\
	&\le U(\mu_n) - U(g_{1/\pi}) - V(\nu_n) + V(g_{1/\pi}) + \delta + \frac{1}{n}\\
	&\le U(\delta_0) - U(g_{1/\pi}) - V(\delta_0) + V(g_{1/\pi}) + \delta + \frac{1}{n}\\
	&\qquad + \Lip(U,d_1) \int_\Rd |x| \mu_n(dx) + \Lip(V,d_1) \int_\Rd |x| \nu_n(dx).
\end{align*}
Upon rearranging terms and using the Cauchy-Schwarz inequality, we find that $\mcl M_2(\mu_n)$ and $\mcl M_2(\nu_n)$ are bounded independently of $n$, and it follows that there exist subsequences $(\mu_{n_k})_{k = 1}^\oo$ and $(\nu_{n_k})_{k =1}^\oo$ and $\bar \mu,\bar \nu \in \mcl P_2(\R^d)$ such that, as $k \to \oo$, $\mu_{n_k} \to \bar \mu$ and $\nu_{n_k} \to \bar \nu$ with respect to the $d_1$ metric, and thus, in particular, weakly in $\mcl P(\R^d)$. In view of the weak lower-semicontinuity of $(\mu, \nu)\mapsto \frac{1}{2 \epsilon }d^2_2(\mu,\nu) + \delta \mathcal{E}^*(\mu) + \delta \mathcal{E}^*(\nu)$ implied by Lemma \ref{entropy:coercive}, as well as the $d_1$-Lipschitz regularity of $U$ and $V$, we obtain
\[
	\Phi(\bar \mu, \bar \nu) \ge \limsup_{k \to \oo} \Phi(\mu_n, \nu_n) \ge \sup \Phi,
\]
and so $(\bar \mu, \bar \nu)$ is a maximum point of $\Phi$ with $|\mcl E(\bar \mu)| + |\mcl E(\bar \nu)| < +\oo$. In particular, $\bar \mu$ and $\bar \nu$ are absolutely continuous with respect to Lebesgue measure. Moreover, since a minimum is achieved by
\[
	\mu \mapsto \mcl E(\mu) + \frac{1}{\delta} \left(- U(\mu) + \frac{1}{2\epsilon} d^2_2(\mu, \bar \nu) + \pi \delta \int_{\R^d} |x|^2 \mu(dx) \right),  
\]
at $\bar{\mu}$, it follows from Lemmas \ref{lem:semicty} and \ref{L:entropyreg} that $\bar \mu^{-1} \in L^\oo_\loc$, and similarly $\bar \nu^{-1} \in L^\oo_\loc$.
 
Given $\mu \in \mathcal{P}_2(\R^d)$, in view of the fact that $\bigl( T_{\bar{\mu}}^{\mu}, T_{\bar{\mu}}^{\bar{\nu}} \bigr)_\# \bar{\mu}$ is an admissible transport plan between $\mu$ and $\bar{\nu}$, we have
\begin{align*}
d_2^2 \bigl( \mu, \bar{\nu} \bigr) &\leq  \int_{\R^d} \bigl| T_{\bar{\mu}}^{\mu}(x) - T_{\bar{\mu}}^{\bar{\nu}}(x) \bigr|^2 d \bar{\mu}(x) \\
& = d_2^2 \bigl( \bar{\mu}, \bar{\nu} \bigr) + 2\int_{\R^d} \bigl( T_{\bar{\mu}}^{\mu}(x)-x \bigr) \cdot \bigl(x- T_{\bar{\mu}}^{\bar{\nu}}(x) \bigr) d \bar{\mu}(x) + d^2_2 \bigl( \bar{\mu},\mu \bigr)
\end{align*}
from which we deduce, thanks to Lemma \ref{L:tangent:crit}, that $p \in \del^+ \bigl(U - \delta \mathcal{E} \bigr) (\oline \mu)$ and $q \in \del^- \bigl(V + \delta \mathcal{E} \bigr) (\oline \nu)$. The $d_1$-Lipschitz continuity of $U$ and $V$ and Lemma \ref{lem:boundedp} together imply the desired $L^\oo$ bounds in \eqref{boundedp}. Since
\[
	p = \frac{1}{\epsilon}( i_d - T_{\oline \mu}^{\oline \nu} ) + 2 \pi \delta i_d \in L^2(\bar{\mu}),
\]
we also have $\nabla  \log \oline \mu \in L^2(\oline \mu)$. Therefore $\mathcal{I}(\oline \mu) < +\infty$, and similarly $\mcl I(\oline \nu) < +\infty$.

Because the transport map $T_{\oline \mu}^{\oline \nu}$ is locally bounded, this also implies that $\nabla \log \oline \mu$ and $\nabla \log \oline \nu$ are locally bounded, and so, since $\oline \mu^{-1}, \oline\nu^{-1} \in L^\oo_\loc$, we conclude that $\bar \mu,\bar \nu \in W^{1,\oo}_\loc$. Finally, the regularity of the transport maps $T_{\bar \mu}^{\bar \nu}$ and $T_{\bar \nu}^{\bar \mu}$ follows from Lemma \ref{L:bootstrapT}.
\end{proof}

\section{Uniqueness and stability}

\label{sec:UniquenessSection}




\subsection{The proof of the comparison principle}

In this section we give the proof of Theorem \ref{UniquenessTheorem3juillet2023}.


\begin{proof}[Proof of Theorem \ref{UniquenessTheorem3juillet2023}]

	By adding a constant to $U$ or $V$, we may assume without loss of generality that $U(T,\cdot) \le V(T,\cdot)$.

	{\it Step 1: Basic penalizations.} For $\gamma > 0$, define
	\[
		U^\gamma(t,\mu) := U(t,\mu) - \gamma(T-t) - \gamma\left( \frac{1}{t} - \frac{1}{T}\right)
	\]
	It is straightforward to check that 
	\begin{equation}\label{penU}
		\left\{
		\begin{split}
		&U^\gamma(T,\cdot) = U(T,\cdot), \quad U^\gamma \le U, \quad \lim_{t \to 0^+} U^\gamma = -\oo, \quad \text{and}\\
		&\text{$U^\gamma$ is a subsolution of \eqref{eq:HJB} with right-hand side $-\gamma$.}
		\end{split}
		\right.
	\end{equation}
	
	Assume for the sake of contradiction that 
	\begin{equation}\label{firstmax}
		[0,T] \ni t \mapsto \sup_{\mu \in \mcl P(\Rd)} \left( U^\gamma(t,\mu) - V(t,\mu) \right)
	\end{equation}
	attains a strictly positive maximum value $\omega > 0$ at some $t_* \in [0,T]$. It follows that $t_* < T$, and, in view of the penalizing property \eqref{penU}, it must also hold that $t_* > 0$.
		
	{\it Step 2: Entropy and diagonal penalizations.} Now let $\epsilon,\delta ,\eta > 0$ be fixed. We then invoke Proposition \ref{P:JKO}, which implies that there exist $(\bar{s} , \bar{t}, \bar{\mu}, \bar{\nu})$ maximizing the function
	\begin{align*}
		\Phi_{\delta,\epsilon}(s,t,\mu,\nu) := U^{\gamma}(s,\mu) - V(t,\nu) - \frac{1}{2\epsilon} d_2^2(\mu,\nu) - \delta  \mathcal{E}^*(\mu) -\delta \mathcal{E}^*(\nu) 
		 - \frac{1}{2\eta}|t-s|^2 
	\end{align*}
	over $[0,T] \times [0,T]  \times \mathcal{P}_2(\R^d) \times \mathcal{P}_2(\R^d)$. Since $\mcl E^* \ge 0$,
	\[
		\delta \to M_{\delta,\epsilon} := \sup_{s,t,\mu,\nu \in [0,T]^2 \times \mcl P_2^2} \Phi_{\delta,\epsilon}(s,t,\mu,\nu)
	\]
	is nonincreasing. It follows that, as $\delta \to 0$ and for fixed $\epsilon > 0$, there exists a limit
	\[
		M_\epsilon := \lim_{\delta \to 0} M_{\delta,\epsilon}.
	\]
	From the inequality
	\begin{align*}
		M_{\delta,\epsilon} = \Phi_{\delta,\epsilon}(\bar s, \bar t, \bar \mu, \bar \nu) \le M_{\delta/2,\epsilon} - \frac{ \delta }{2} \mcl E^*(\bar \mu) - \frac{\delta}{2} \mcl E^*(\bar \nu)
	\end{align*}
	and the fact that $M_{\delta/2,\epsilon} - M_{\delta,\epsilon} \xrightarrow{\delta \to 0} 0$, we conclude that, for fixed $\epsilon > 0$, 
	\[
		\lim_{\delta \to 0} \delta \left( \mcl E^*(\bar \mu) + \mcl E^*(\bar \nu) \right) = 0.
	\]
	With another application of \eqref{entropy:coercive}, this also implies, for fixed $\epsilon > 0$,
	\begin{equation}\label{moments:small}
		\lim_{\delta \to 0}\delta \left(  \int |x|^2 d\bar\mu(x) +  \int |y|^2 d\bar\nu(x)\right)=0.
	\end{equation}
Rearranging terms in the inequality $\Phi(\bar s, \bar t, \bar \nu, \bar \nu) + \Phi(\bar s, \bar t, \bar \mu, \bar \mu)  \le 2\Phi(\bar s, \bar t, \bar \mu, \bar \nu)$ leads to
	\[
		\frac{1}{2\epsilon} d^2_2(\bar \mu, \bar \mu) \le U(\bar s, \bar \mu) - U(\bar s, \bar \nu) + V(\bar t, \bar \mu) - V(\bar t, \bar \nu)
		\le (L_U + L_V)d_1(\bar \mu, \bar \nu) \le (L_U + L_V)d_2(\bar \mu, \bar \nu),
	\]
	from which we conclude
	\begin{equation}\label{diagonal:small}
		d_2(\bar{\mu}, \bar{\nu}) \leq 2(L_U + L_V) \epsilon. 
	\end{equation}
	It follows that
	\[
		\lim_{\eta \to 0} \lim_{\epsilon \to 0} \lim_{\delta \to 0} \Phi_{\delta,\epsilon}(\bar s, \bar t, \bar \mu, \bar \nu) = \omega,
	\]
	and so, for all sufficiently small $\epsilon, \delta,\eta > 0$,
	\[
		\Phi_{\delta,\epsilon}(\bar s, \bar t, \bar \mu, \bar \nu) \ge \frac{\omega}{2}.
	\]
	Moreover, as $\epsilon \to 0$, $\delta \to 0$, and then $\eta \to 0$, any limit point of $\bar s$ and $\bar t$ must agree, and be a maximizing point of the original function \eqref{firstmax}. We deduce that, for $\epsilon, \delta,\eta$ sufficiently small, we have $0 < \bar s < T$ and $0 < \bar t < T$.
	
	{\it Step 3: Applying the definition.} By Proposition \ref{P:JKO}, $\mcl I(\bar \mu) + \mcl I(\bar \nu) < \oo$, and we may therefore apply Definition \ref{D:soln} to obtain the inequalities
	\begin{align*}
		-\frac{\bar s- \bar t}{\eta} &+ \int_{\R^d} H \Bigl(x,\frac{1}{\epsilon}(x-T_{\bar{\mu}}^{\bar{\nu}}(x)) +2\pi \delta x, \bar{\mu} \Bigr)d\bar{\mu}(x) \\
		&+\int_{\R^d} \left( \div_x a(\bar\mu,x) + a(\bar\mu,x)\nabla \log \bar{\mu}(x)\right) \cdot \Bigl( \frac{1}{\epsilon}(x-T_{\bar{\mu}}^{\bar{\nu}}(x)) +2\pi \delta x \Bigr) d\bar \mu(x) \leq -\gamma + C\delta,
	\end{align*}
	and
	\begin{align*}
		-\frac{\bar s-\bar t}{\eta} &+ \int_{\R^d} H \Bigl(y, \frac{1}{\epsilon}(T_{\bar{\nu}}^{\bar{\mu}}(y)-y)  - 2\pi \delta y, \bar{\nu} \Bigr)d\bar{\nu}(y) \\
		&+ \int_{\R^d} \left( \div_y a(\bar \nu,y) + a(\bar \nu,y)\nabla \log \bar{\nu}(y) \right)\cdot \Bigl( \frac{1}{\epsilon}(T_{\bar{\nu}}^{\bar{\mu}}(y)-y) - 2\pi \delta y \Bigr) d\bar \nu(y) \geq - C\delta.
	\end{align*}
Here $C$ is the maximum of the two constants for $U$ and $V$ from Definition \ref{D:soln}, and, in what follows, the constant $C$ may change from line to line, and is independent of $\epsilon$, $\delta$, and $\gamma$.
	
	Denote by $\bar{\pi} = (\Id, T_{\bar\mu}^{\bar\nu})_\# \bar \mu$ the optimal transport plan between $\bar{\mu}$ and $\bar{\nu}$. Subtracting the two inequalities, we obtain
	\[
		\gamma \le 2C \delta + I_1 + I_2 + I_3,	
	\]
	where\footnote{Here we are using the fact that $\bar\mu$ and $\bar\nu$ are positive and Lipschitz on $\Rd$, so that the pointwise identities $\nabla \log \bar \mu(x) \bar \mu(x) = \nabla \bar \mu(x)$ and $\nabla \log \bar \nu(x) \bar \nu(x) = \nabla \bar \nu(x)$ are justified.}
	\begin{align*}
		I_1 &:= \int_{\R^d \times \R^d} \Bigl[ H \Bigl( y, \frac{1}{\epsilon}(x-y)  - 2\pi \delta y, \bar{\nu} \Bigr) - H \Bigl( x, \frac{1}{\epsilon}(x-y) +2\pi \delta x, \bar{\mu} \Bigr)\Bigr] d\bar{\pi}(x,y),
	\end{align*}
	\begin{align*}
		I_2 &:=  \frac{1}{\epsilon} \int_{\R^d} \left( \div_x a(\bar\mu,x) + a(\bar\mu,x) \nabla \bar{\mu}(x)\right) \cdot (T_{\bar{\mu}}^{\bar{\nu}}(x) - x)  dx\\
		&- \frac{1}{\epsilon} \int_{\R^d} \left( \div_y a(\bar \nu,y) + a(\bar \nu,y)\nabla \bar{\nu}(y) \right)\cdot (y - T_{\bar{\nu}}^{\bar{\mu}}(y)) dy,
	\end{align*}
	and
	\begin{align*}
		I_3 := -2\pi \delta \int_{\Rd}  \left( \div_x a(\bar \mu, x) + a(\bar \mu,x) \nabla \bar \mu(x) + \div_x a(\bar \nu,x) + a(\bar \nu,x) \nabla \bar \nu(x) \right) \cdot x dx.
	\end{align*}
	In view of the assumptions \eqref{AHexistence 3juillet2023} on $H$, there exists $C > 0$ depending only on the assumptions for $H$, such that, by several applications of the Cauchy-Schwarz inequality,
	\begin{align*}
		|I_1| &\le C \int_{\Rd \times \Rd} \left( 1 + \frac{|x-y|}{\epsilon} + \delta|x| + \delta|y| \right)\left( |x-y| + \delta|x| + \delta|y| + d_2(\bar \mu, \bar \nu) \right)d \bar \pi(x,y) \\
		&\le C\left( d_2(\bar \mu, \bar \nu) + \frac{d_2^2(\bar \mu, \bar \nu)}{\epsilon} + \delta + \delta^2 \int_\Rd |x|^2 d \bar \mu(x) + \delta^2 \int_\Rd |y|^2 d\bar \nu(y) \right)\\
		&\le C\epsilon + C \delta,
	\end{align*}
	where the last line is a consequence of \eqref{moments:small} and \eqref{diagonal:small}.

	Now let $\chi \in C_c^\oo(\Rd)$ be a nonnegative bump function satisfying $\chi = 1$ in $B_1$ and $\chi = 0$ in $\Rd \backslash B_2$, and, for $R > 0$, define $\chi_R(x) = \chi(x/R)$. Set
	\begin{align*}
		I_{2,R} &:=  \frac{1}{\epsilon} \int_{\R^d} \chi_R(x)\left( \div_x a(\bar\mu,x) + a(\bar\mu,x) \nabla \bar{\mu}(x)\right) \cdot (T_{\bar{\mu}}^{\bar{\nu}}(x) - x)  dx\\
		&- \frac{1}{\epsilon} \int_{\R^d}\chi_R(y) \left( \div_y a(\bar \nu,y) + a(\bar \nu,y)\nabla \bar{\nu}(y) \right)\cdot (y - T_{\bar{\nu}}^{\bar{\mu}}(y)) dy.
	\end{align*}
	The dominated convergence theorem implies that $\lim_{R \to \oo} I_{2,R} = I_2$. On the other hand, in view of the compact support of $\chi_R$ and the $C^1$ regularity of $T_{\bar \mu}^{\bar \nu}$, we may integrate by parts to obtain
	\begin{equation}\label{I2Rab}
	\begin{split}
		I_{2,R} 
		&= \frac{1}{\epsilon} \int_\Rd \chi_R(x) \tr\left[ a(\bar \mu,x)(I_d -\nabla T_{\bar{\mu}}^{\bar{\nu}}(x)) \right] d\bar \mu(x) - \frac{1}{\epsilon} \int_\Rd \chi_R(x) \tr\left[ a(\bar \nu, y) (\nabla T_{\bar \nu}^{\bar \mu}(x) - I_d) \right] d \bar \nu(x)\\
		&+ \frac{1}{\epsilon} \int_\Rd \nabla \chi_R(x)a(\bar \mu,x)(x-T_{\bar{\mu}}^{\bar{\nu}}(x))d \bar \mu(x) - \frac{1}{\epsilon} \int_\Rd \nabla \chi_R(y) a(\bar \nu, y) (T_{\bar \nu}^{\bar \mu}(y) - y)d \bar \nu(y).
	\end{split}
	\end{equation}
	By the Cauchy-Schwarz inequality, the terms in the second line of \eqref{I2Rab} can be estimated, for some constant $C > 0$, by
	\begin{equation}\label{I2Rb}
		\begin{split}
			&\frac{1}{\epsilon} \left| \int_\Rd \nabla \chi_R(x)a(\bar \mu,x)(x-T_{\bar{\mu}}^{\bar{\nu}}(x))d \bar \mu(x)\right| +  \frac{1}{\epsilon} \left|\int_\Rd \nabla \chi_R(y) a(\bar \nu, y) (T_{\bar \nu}^{\bar \mu}(y) - y)d \bar \nu(y)\right|\\
			&\le \frac{C\norm{a}_\oo}{R\epsilon} d_2(\bar \mu, \bar \nu).
		\end{split}
	\end{equation}
	We have $\nabla T_{\bar \nu}^{\bar \mu}(y) = \nabla T_{\bar \mu}^{\bar \nu}(x)^{-1}$ for $\bar \pi$-a.e. $(x,y) \in \Rd$, and so the first line of \eqref{I2Rab} can be rewritten as
	\begin{align*}
		&\frac{1}{\epsilon} \int_\Rd \chi_R(x) \tr\left[ a(\bar \mu,x)(I_d -\nabla T_{\bar{\mu}}^{\bar{\nu}}(x)) \right] d\bar \mu(x) - \frac{1}{\epsilon} \int_\Rd \chi_R(x) \tr\left[ a(\bar \nu, y) (\nabla T_{\bar \nu}^{\bar \mu}(x) - x) \right] d \bar \nu(x)\\
		&= \frac{1}{\epsilon} \int_{\Rd \times \Rd} \left( \chi_R(x) \tr\left[a(\bar \mu, x) (I_d - \nabla T_{\bar \mu}^{\bar \nu}(x) )\right] + \chi_R(y)  \tr\left[ a(\bar \nu, y) (I_d - \nabla T_{\bar \mu}^{\bar \nu}(x)^{-1}) \right] \right)d\bar \pi(x,y).
	\end{align*}
	The convexity of $\phi_{\bar \mu}^{\bar \nu}$ and $\phi_{\bar \nu}^{\bar \mu}$ then gives, for all $R$,
	\begin{align*}
		&\chi_R(x) \tr\left[a(\bar \mu, x) (I_d - \nabla T_{\bar \mu}^{\bar \nu}(x) )\right] + \chi_R(y)  \tr\left[ a(\bar \nu, y) (I_d - \nabla T_{\bar \mu}^{\bar \nu}(x)^{-1}) \right]\\
		&\le \tr[a(\bar \mu, x)] + \tr[a(\bar \nu,y)] \in L^1(\Rd \times \Rd, \bar \pi),
	\end{align*}
	and so, by Fatou's lemma,
	\begin{equation}\label{matrixtime}
	\begin{split}
		&\limsup_{R \to \oo} \frac{1}{\epsilon} \int_{\Rd \times \Rd} \left( \chi_R(x) \tr\left[a(\bar \mu, x) (I_d - \nabla T_{\bar \mu}^{\bar \nu}(x) )\right] + \chi_R(y)  \tr\left[ a(\bar \nu, y) (I_d - \nabla T_{\bar \mu}^{\bar \nu}(x)^{-1}) \right] \right)d\bar \pi(x,y)\\
		&\le \frac{1}{\epsilon} \int_{\Rd \times \Rd} \left(  \tr\left[a(\bar \mu, x) (I_d - \nabla T_{\bar \mu}^{\bar \nu}(x) )\right] +  \tr\left[ a(\bar \nu, y) (I_d - \nabla T_{\bar \mu}^{\bar \nu}(x)^{-1}) \right] \right)d\bar \pi(x,y).
	\end{split}
	\end{equation}
	For a positive-definite $d\times d$ symmetric matrix $X$ and two $d \times d$ matrices $\alpha$ and $\beta$, rearranging the inequality
	\[
		\tr[ (X^{1/2} \alpha - X^{-1/2} \beta)(X^{1/2} \alpha - X^{-1/2} \beta)^* ] \ge 0
	\]
	yields
	\begin{equation}\label{linalg}
		\tr[ X\alpha\alpha^*] + \tr[ X^{-1} \beta \beta^*] \ge 2 \tr[ \alpha \beta^*].
	\end{equation}
	Applying \eqref{linalg} with $X = \nabla T_{\bar \mu}^{\bar \nu}(x)$, $X^{-1} = \nabla T_{\bar{\nu}}^{\bar{\mu}}(y)$, $\alpha = \sigma(\bar \mu,x)$, and $\beta = \sigma(\bar \nu,y)$, the last line of \eqref{matrixtime} can be further developed to, for some $C > 0$ depending on the constants in \eqref{A:a},
	\begin{align*}
		&\frac{1}{\epsilon}  \int_{\R^d \times \R^d} \tr \bigl(   a(\bar \mu,x) +  a(\bar \nu,y)  - a(\bar \mu,x)\nabla T_{\bar{\mu}}^{\bar{\nu}}(x) - a(\bar \nu,y) \nabla T_{\bar{\nu}}^{\bar{\mu}}(y) \bigr) d\bar{\pi}(x,y)  \\
		&\le \frac{1}{\epsilon} \int_{\Rd \times \Rd}  \tr \left( \sigma(\bar \mu, x)\sigma(\bar \mu, x)^* + \sigma(\bar \nu, y)\sigma(\bar\nu,y)^* - 2\sigma(\bar \mu, x)\sigma(\bar \nu, y)^*  \right) d \bar \pi(x,y) \\
		&=  \frac{1}{\epsilon} \int_{\Rd \times \Rd} \tr \left[ (\sigma(\bar \mu,x) - \sigma(\bar \nu, y) )(\sigma(\bar \mu,x) - \sigma(\bar \nu, y) )^* \right]d \bar \pi(x,y) \le  \frac{C}{\epsilon} d_2^2(\bar \mu, \bar \nu) .
	\end{align*}
	Combining the above lines with \eqref{I2Rb}, we conclude that
	\[
		I_2 = \lim_{R \to \oo} I_{2,R} \le \frac{C}{\epsilon} d_2^2(\bar \mu, \bar \nu) \le C\epsilon.
	\]
	Finally, integrating by parts gives, for some constant $C > 0$,
	\[
		I_3 = 2\pi \delta \int_\Rd \tr[a(\bar \mu,x)]d \bar \mu(x) + 2\pi \delta \int_\Rd \tr[ a(\bar \nu,y)]d \bar \nu(y) \le C \delta.
	\]
	Combining the estimates for $I_1$, $I_2$, and $I_3$, we obtain
	\begin{equation}\label{last:inequality}
		\gamma \le C\epsilon + C\delta.
	\end{equation}
	Sending first $\delta \to 0$ and then $\epsilon \to 0$ gives the desired contradiction.
\end{proof}


The proof of the comparison principle can be easily modified to obtain the following stability estimate.

\begin{thm}\label{T:stability}
	Assume $(H^1,H^2)$ satisfy \eqref{AHexistence 3juillet2023} and $\mathcal{G}_1, \mathcal{G}_2$ are two Lipschitz continuous maps w.r.t. $d_1$, and let $U^1$ and $U^2$ be corresponding sub and supersolutions satisfying the conditions of Theorem \ref{UniquenessTheorem3juillet2023}, whose $d_1$-Lipschitz norms are bounded by $L > 0$. Then, for all $(t, \mu) \in [0,T] \times \mcl P_2(\R^d)$,
	\begin{align*}
		 U^1(t,\mu) - U^2(t,\mu) 
		&\le \sup_{\mu \in \mcl P_2(\R^d)} \left( \mathcal{G}_1(\mu) - \mathcal{G}_2(\mu) \right)_+\\
		&+ (T-t) (1 + 8L^2) \sup_{(x,p,\mu) \in \Rd \times \Rd \times \mcl P_2(\R^d)} \frac{\left[ H_2(x,p,\mu) - H_1(x,p,\mu) \right]_+}{1 + |p|^2}.
	\end{align*}
\end{thm}

\begin{proof}
	The proof follows almost exactly as for Theorem \ref{UniquenessTheorem3juillet2023} with $U = U^1$ and $V = U^2$. The main difference is the presence of the error term involving the difference between $H_1$ and $H_2$. Setting
	\[
		M := \sup_{(x,p,\mu) \in \Rd \times \Rd \times \mcl P_2(\R^d)} \frac{\left[ H_2(x,p,\mu) - H_1(x,p,\mu) \right]_+}{1 + |p|^2},
	\]
	the inequality \eqref{last:inequality} then reads as
	\begin{align*}
		\gamma &\le C\delta + C\epsilon + \int_{\Rd} \left[ H_2\left( x, \frac{x- T_{\bar \mu}^{\bar \nu}(x)}{\epsilon} ,\bar \mu \right) - H_1 \left( x, \frac{x-T_{\bar \mu}^{\bar \nu}(x)}{\epsilon}, \bar \mu\right) \right]_+ d\bar \mu(x)\\
		&\le C \delta + C\epsilon + M \int_{\Rd}\left(1 + \frac{1}{\epsilon^2} |x - T_{\bar \mu}^{\bar \nu}(x)|^2\right) d \bar \mu(x) \\
		&= C \delta + C\epsilon + M + \frac{M}{\epsilon^2} d_2^2(\bar \mu, \bar \nu) \le C \delta + C\epsilon + M + 8ML^2,
	\end{align*}
	where the last inequality is a consequence of \eqref{diagonal:small}. The proof is finished upon sending $\delta \to 0$ and $\epsilon \to 0$.
\end{proof}

\subsection{A posteriori $d_2$-Lipschitz estimate}

Our proof of the comparison principle requires a priori $d_1$-Lipschitz regularity for $U$ and $V$, which is stronger than $d_2$-Lipschitz regularity. On the other hand, we now demonstrate that the $d_2$-Lipschitz constant of a solution can be controlled by the $d_2$-Lipschitz constant of the terminal datum $G$, independently of the $d_1$-Lipschitz constant. 
	
\begin{prop}\label{P:d2Lip}
	There exists a constant $C > 0$ depending only on $T$ and the constants in the assumptions \eqref{A:a} and \eqref{AHexistence 3juillet2023} such that, if $U$ is a $d_1$-Lipschitz viscosity solution of \eqref{eq:HJB}, then
	\[
		\sup_{t \in [0,T]} \Lip(U(t,\cdot); d_2) \le C (1 + \Lip(G; d_2)).
	\]
\end{prop}

\begin{rmk}
	As will be clear from the proof, the constant $C$ above is also independent of the constant in front of $\delta$ in Definition \ref{D:soln}. Recall from Remark \ref{rmk:elliptic} that the constant $C$ should be expected to depend on the ellipticity constant of $a$, that is, the parameter $\kappa$ in \eqref{A:aelliptic}. It follows that the above Lipschitz estimate is itself independent of the ellipticity constant. 
\end{rmk}

\begin{proof}[Proof of \ref{P:d2Lip}]
	Fix $0 < \tau < T$, and set 
	\[
		L := \sup_{t \in [T-\tau,T]} \Lip(U(t,\cdot);d_2),
	\]
	which is finite in view of the $d_1$-Lipschitz regularity of $U$. For $\alpha > 0$, set $\tilde U(t,\mu) := U(t,\mu) - \alpha (t- T + \tau)^{-1}$, and, for $\epsilon, \gamma > 0$, suppose that
	\[
		t \mapsto \sup_{\mu, \nu \in \mcl P_2(\R^d)} \left\{ \tilde U(t,\mu) - U(t,\nu) - \frac{1}{2\epsilon} d^2_2(\mu,\nu) \right\} - \gamma (T-t)
	\]
	attains a maximum in $[T-\tau,T]$ at some $t_0 < T$, which, because $\alpha > 0$, must also satisfy $t_0 > T - \tau$. Then, arguing just as in the proof of Theorem \ref{UniquenessTheorem3juillet2023}, it follows that, for $\delta, \eta > 0$, the map
	\[
		(s,t,\mu,\nu) \mapsto \tilde U(s,\mu) - U(t,\nu) - \frac{1}{2\epsilon} d_2^2(\mu,\nu) - \delta \mathcal{E}^*(\mu) - \delta \mathcal{E}^*(\nu) - \frac{1}{2\eta}|t-s|^2 - \gamma (T-t) - \frac{1}{2} |s-t_0|^2
	\]
	attains a maximum $(\bar s, \bar t, \bar \mu, \bar \nu) \in [T-\tau,T]^2 \times \mcl P_2(\R^d)^2$ such that $\bar \mu$ and $\bar \nu$ satisfy the conclusions of Proposition \ref{P:JKO}. Furthermore, \eqref{moments:small} holds uniformly in $\eta$, $\lim_{\eta \to 0} \lim_{\delta \to 0} (\bar s, \bar t) = (t_0,t_0)$ (and therefore $\bar s$ and $\bar t$ lie strictly inside $(T - \tau, T)$ for sufficiently small $\delta$ and $\eta$), and
	\begin{equation}\label{L:d2}
		d_2(\bar \mu, \bar \nu) \le L\epsilon.
	\end{equation}
	We then argue similarly as in the proof of Theorem \ref{T:comparison}. Keeping track of the dependence of the estimates, using \eqref{L:d2}, we have the same error estimates, for some constant as in the statement in the theorem,
	\[
		|I_1| \le C(1 + L^2)\epsilon + C \delta \quad \text{and} \quad |I_2| \le CL^2 \epsilon.
	\]
	Denoting by $C_0$ the constant\footnote{Specifically, $C_0$ is the maximum of the constants for $U$ and for $V$.} from Definition \ref{D:soln}, we arrive at the inequality
	\[
		\gamma - (\bar s - t_0) \le C_0 \delta + C \delta + C(1 + L^2)\epsilon.
	\]
	Sending $\delta \to 0$ and $\eta \to 0$ gives a contradiction if $\gamma > C(1 + L^2)\epsilon$. Since $\alpha$ was arbitrary, we thus have, for all $t \in [T-\tau,T]$ and $\mu,\nu \in \mcl P_2$,
	\begin{align*}
		U(t,\mu) - U(t,\nu) - \frac{1}{2\epsilon} d^2_2(\mu,\nu) &\le C\tau(1 + L^2)\epsilon + \sup_{\mu,\nu \in \mcl P_2} \left\{ G(\mu) - G(\nu) - \frac{1}{2\epsilon}d_2^2(\mu,\nu) \right\}\\
		&\le C\tau(1 + L^2)\epsilon + \frac{1}{2} \Lip(G;d_2)^2 \epsilon.
	\end{align*}
	Optimizing in $\epsilon$ gives, for a possibly different $C$,
	\[
		U(t,\mu) - U(t,\nu) \le C\left( \tau^{1/2}(1 + L) + \Lip(G;d_2) \right) d_2(\mu,\nu),
	\]
	and therefore
	\[
		L \le C\left( \tau^{1/2}(1 + L) + \Lip(G;d_2) \right) .
	\]
	If $\tau$ is sufficiently small, depending only on $C$, we arrive at
	\[
		\Lip(U(t,\cdot);d_2) \le C(1 +  \Lip(G;d_2)) \quad \text{for } t \in [T-\tau,T].
	\]
	Because $\tau$ is independent of $G$, this result can then be iterated on intervals of the form $[T-2\tau,T-\tau]$, $[T-3\tau,T-2\tau]$, etc., yielding the Lipschitz estimate on the full interval $[0,T]$ with a possibly different constant $C$ depending additionally on $T$.
\end{proof}


\section{Existence of solutions and optimal control of the Fokker-Planck equation}

\label{sec:ExistenceSection}

The goal of this section is to prove the existence result Theorem \ref{thm:ExistenceThm}. We assume that $H$ takes the form
$$H(x,p,\mu) = H_1(x,p) - \mathcal{F}(\mu),$$ 
with $H_1$ convex in the second variable. In this case, the Hamilton-Jacobi equation \eqref{eq:HJB} becomes
\begin{equation}
\left \{
    \begin{array}{ll}
  \displaystyle       -\partial_t U(t,\mu) - \int_{\R^d}  \div_x \partial_{\mu} U(t,\mu,x) d\mu(x)  \\
 \displaystyle     \hspace{100pt} + \int_{\R^d} H_1 \bigl(x, \partial_{\mu}U(t,\mu,x) \bigr) d\mu(x) = \mathcal{F}(\mu)  &  \mbox{ in } [0,T] \times \mathcal{P}_2(\R^d)\\
    U(T,\mu) = \mathcal{G}(\mu)     & \mbox{in } \mathcal{P}_2(\R^d). 
    \end{array}
\right.
\label{HJBExistence}
\end{equation}
We assume that $\mathcal{F}, \mathcal{G}$ and $H_1$ satisfy \eqref{AHallinone} and \eqref{AssumptionsMeanFieldCosts} from Section \ref{sec:Defofviscositysolandmainresults}, and we prove that the unique solution of \eqref{HJBExistence} is given by the value function $U : [0,T] \times \mathcal{P}_2(\R^d) \rightarrow \R$ of the optimal control problem. More precisely, we set $L(x,q) := \sup_{p \in \R^d} \bigl \{ -p \cdot q - H_1(x,p) \bigr \}$, and, for $(t_0,m_0) \in [0,T] \times \mathcal{P}_2(\R^d)$, 
\begin{equation}\label{def:Uexistence}
	U(t_0,m_0) = \inf_{(\mu,\alpha)} \left\{ \int_{t_0}^T \int_{\R^d} L \bigl(x,\alpha_t(x) \bigr)d\mu_t(x)dt +\int_{t_0}^T \mathcal{F}(\mu_t)dt + \mathcal{G}(\mu_T)\right\},
\end{equation}
where the infimum is taken over couples $(\mu,\alpha)$ with $\mu \in \mathcal{C}([t_0,T] ,\mathcal{P}_2(\R^d)), \alpha\in L^2 ([t_0,T] \times \R^d, dt\otimes \mu_t; \R^d) $ satisfying, in the sense of distributions, the Fokker-Planck equation
\begin{equation}\label{FP}
\left \{
\begin{array}{ll}
\partial_t \mu + \div((\alpha_t(x)\mu)-\Delta \mu = 0 &\mbox{ in } (t_0,T) \times \R^d,  \\
\mu(t_0)=\mu_0.
\end{array}
\right.
\end{equation}

We first establish the boundedness and regularity of the value function. As a byproduct of the proof, the infimum in the definition of $U$ can be restricted to uniformly bounded controls, a fact which is crucial in verifying the sub and supersolution properties of $U$. We omit the proof, as it follows exactly as in \cite{DDJ_23} (see also \cite{CDJS_23}).

\begin{prop}
	Assume \eqref{AHallinone} and \eqref{AssumptionsMeanFieldCosts}. Then $U$ is bounded, continuous and Lipschitz with respect to $d_1$ uniformly in time. Moreover, there exists $C>0$ independent of $(t_0,\mu_0) \in [0,T] \times \mathcal{P}_2(\R^d)$ such that the infimum in the definition of $U(t_0,\mu_0)$ can be restricted to controls satisfying  $\displaystyle \|\alpha \|_{\infty} \leq C.$
\label{prop: regularity of the value function}
\end{prop}

Since we will be working with bounded controls in what follows, we will be able to take advantage of the following continuity result. We provide its short proof for the sake of completeness.

\begin{lem}\label{lem:L1cts}
	Assume $\alpha \in L^\oo([t_0,T] \times \Rd)$. Then the solution $\mu$ of \eqref{FP} belongs to $C([t_0,T], L^1(\Rd))$, with a modulus of continuity depending only on $\norm{\alpha}_\oo$ and $\mu_0$.
\end{lem} 


\begin{proof}
	Take $t_2,t_1 \in [t_0,T]$ with $t_2 \geq t_1$. The solution $\mu_t$ can be expressed with Duhamel's formula by
	\[
		\mu_{t_2} = p_{t_2 - t_1} * \mu_{t_1} - \int_{t_1}^{t_2} (\nabla p_{t-s}) * (\alpha_s \mu_s)ds,
	\]
where 
\begin{equation}\label{heat}
	p_t(x) := \frac{1}{(4\pi t)^{d/2}} \exp\left( - \frac{|x|^2}{4t} \right), \quad (t,x) \in (0,\oo) \times \Rd.
\end{equation}
is the heat kernel.
%
Subtracting $\mu_{t_1}$ from both sides and taking the $L^1$-norm yields, using Young's convolution inequality, for some dimensional constant $c_d > 0$,
%
%
%
	\begin{align*}
		\norm{\mu_{t_2} - \mu_{t_1}}_{L^1(\Rd)} & \leq \norm{ p_{t_2-t_1} * \mu_{t_1} - \mu_{t_1}}_{L^1} +  \int_{t_1}^{t_2} \norm{\nabla p_{t-s}}_{L^1(\R^d)} \norm{ \alpha_s \mu_s}_{L^1(\R^d)}ds \\
 & \le \norm{ p_{t_2-t_1} * \mu_0 - \mu_0}_{L^1} + \norm{\alpha}_{\infty} \int_{t_1}^{t_2} \norm{ \nabla p_{t-s}}_{L^1(\Rd)}ds\\
		&\le \norm{ p_{t_2-t_1} * \mu_0 - \mu_0}_{L^1} + c_d \norm{\alpha}_{\infty} \int_{t_1}^{t_2}(t-s)^{-1/2}ds\\
		&=\norm{ p_{t_2-t_1} * \mu_0 - \mu_0}_{L^1} + 2c_d \norm{\alpha}_{\infty} (t_2-t_1)^{1/2}. 
	\end{align*}
The proof is finished because $\lim_{t_2 \to t_1} \norm{ p_{t_2-t_1} * \mu_0 - \mu_0}_{L^1} = 0$, at a rate depending only on the $L^1$-modulus of continuity of $\mu_0$.
\end{proof}

When $(\mu_t)_{t\geq t_0}$ is a solution to \eqref{FP}, it is elementary to check, using the stochastic characteristics, that
$$d_2 \bigl( \mu_t, \mu_0 \bigr) \leq \norm{\alpha}_{\infty}|t-t_0| + \sqrt{2} \sqrt{|t-t_0|}, \quad \forall t \geq t_0,$$
which, for arbitrary $\mu_0$, is sharp (if, for instance, $\mu_0$ is a Dirac $\delta$-distribution and $\alpha=0$). In checking that the value function is a viscosity solution, it will be crucial to upgrade the modulus of $|t-t_0|$ to be linear, which turns out to be the case under the additional assumption that the Fischer information of $\mu_0$ is finite. 

\begin{prop}\label{prop:FPresults}
	Fix $(t_0,\mu_0) \in [0,T] \times \mathcal{P}_{2,ac}(\R^d)$ with $\mcl I(\mu_0) < +\oo$, and assume $\mu \in \mathcal{C}([t_0,T] ,\mathcal{P}_2(\R^d))$ and $\alpha\in L^{\infty} ([t_0,T] \times \R^d,\R^d)$ satisfy \eqref{FP} in the sense of distributions. Then the following hold:
\begin{enumerate}[(a)]
\item\label{prop:TimeLipschitzEstimate} For all $t \in [t_0,T]$, 
$$d_2 \bigl (\mu_t,\mu_0 \bigr) \leq \bigl( \norm{\alpha}_{\infty} +  \mathcal{I}(\mu_0)^{1/2} \bigr) |t-t_0|.$$
\item\label{prop:IPP} For $h > 0$, let $T_{\mu_0}^{\mu_{t_0+h}}$ denote the optimal transport map between $\mu_0$ and $\mu_{t_0+h}$. Then there exists $C > 0$ depending only on $\mcl I(\mu_0)$ and $\norm{\alpha}_\oo$ such that, for any $\varphi \in \mathcal{C}^{\infty}_c(\R^d)$, 
\[
	\left| \int_{\R^d} \nabla \varphi(x) \cdot  \left(T_{\mu_0}^{\mu_{t_0+h}}(x) - x \right) d\mu_{0}(x) - \int_{t_0}^{t_0 + h} \int_{\R^d} \bigl[  \alpha_{t}(x)\cdot \nabla \varphi(x) +\Delta \varphi(x) \bigr] d\mu_0(x) \right| \le \norm{\nabla^2 \varphi}_\oo  h^2.
\]
In particular, if $t \mapsto \alpha_t$ is continuous at $t_0^+$ in $L^2(\mu_0)$, then 
\[
	\lim_{h \rightarrow 0^+}\int_{\R^d} \nabla \varphi(x) \cdot \frac{T_{\mu_0}^{\mu_{t_0+h}}(x) - x}{h} d\mu_{0}(x) = \int_{\R^d} \bigl[  \alpha_{t_0}(x)\cdot \nabla \varphi(x) +\Delta \varphi(x) \bigr] d\mu_0(x).
\]
\end{enumerate}
\end{prop}

\begin{proof}
\eqref{prop:TimeLipschitzEstimate} For a fixed probability space $\bigl( \Omega, \mathcal{F}, \mathbb{P} \bigr)$ supporting a $d$-dimensional Brownian motion $(B_t)_{t \geq 0}$, let $X_{t_0}: \Omega \to \RR^d$ be a random variable with law $\mu_{0}$, independent of $B$, and denote by $(X_t)_{t \geq t_0}$ and $(Y_t)_{t \geq t_0}$ the solutions of respectively
$$X_t = X_{t_0} + \int_{t_0}^t \alpha(s,X_s)ds + \sqrt{2} \bigl(B_t- B_{t_0}), \quad t \geq t_0$$
and
$$Y_t = X_{t_0} + \sqrt{2} \bigl(B_t- B_{t_0}), \quad t \geq t_0.$$
Then, for $t \geq t_0$, $\mu_t$ is equal to the law of $X_t$, while the law $\nu_t$ of $Y_t$ solves
\[
	\partial_t \nu_t - \Delta \nu_t =0, \quad \mbox{in } (t_0,T) \times \R^d, \quad \nu_{t_0} = \mu_0.
\]
It follows that
\begin{equation}\label{estimatewithdrift3juillet}
	d_2^2 \bigl(\mu_t, \nu_t \bigr) \leq \E \bigl[  \bigl| X_t - Y_t \bigr|^2 \bigr] 
		\leq \E \Bigl[ \Bigl| \int_{t_0}^t \alpha(s,X_s)ds \Bigr|^2 \Bigr]
		\leq |t-t_0|^2 \norm{\alpha}^2_{\infty}.
\end{equation}

Writing $(P_t)_{t \geq 0}$ for the heat kernel (associated to $-\Delta$) we can rewrite, for all $t \geq t_0$,
$$ \mathcal{I}(\nu_t) = \int_{\R^d} \frac{|\nabla \nu_t(x)|^2}{\nu_t(x)}dx = \int_{\R^d} \frac{\bigl| \nabla P_{t-t_0} \mu_0(x) \bigr|^2}{P_{t-t_0}\mu_0(x)}dx  = \int_{\R^d} \frac{\bigl| P_{t-t_0} \nabla \mu_0(x) \bigr|^2}{P_{t-t_0}\mu_0(x)}dx.$$
We now use Jensen's inequality and the convexity of the map $\R^d \times \R^+ \ni (a,b) \mapsto |a|^2/b$ (set to be $+\infty$ if $b=0$), to estimate, for all $t \geq t_0$,
\begin{equation} 
\mathcal{I}(\mu_t) \leq \int_{\R^d} P_{t-t_0} \frac{ \bigl|\nabla \mu_0(.) \bigr|^2}{\mu_0(.)}(x)dx = \mathcal{I}(\mu_0),
\label{JensenforFischerinformation}
\end{equation}
where we used the conservation of mass along the heat flow for the last equality. 
The heat equation can be rewritten as a continuity equation driven by $-\nabla \log \nu_t$, and so, by Proposition \ref{P:transportplan}\eqref{BenBre},
$$d_2^2 \bigl( \nu_t, \mu_0 \bigr) \leq |t-t_0| \int_{t_0}^t \norm{ \nabla \log \nu_s}^2_{L^2(\nu_s)} ds = |t-t_0| \int_{t_0}^t \mathcal{I}(\nu_s)ds.$$
Estimate \eqref{JensenforFischerinformation} therefore gives $d_2^2 \bigl( \nu_t, \mu_0 \bigr) \leq |t-t_0|^2 \mathcal{I}(\mu_0)$. Combining this estimate with \eqref{estimatewithdrift3juillet} then yields the desired estimate for $d_2(\mu_t,\mu_0)$.

\eqref{prop:IPP} Arguing just as for \eqref{tobeusedlaterinSection5} and \eqref{tobeusedlateraswell} in the proof of Lemma \ref{lem:boundedp}, we have the inequality
\begin{align*}
	&\left| \int_\Rd \nabla \varphi(x) \cdot \left(T_{\mu_0}^{\mu_0 + h}(x) - x \right)d\mu_0(x)
	- \int_{t_0}^{t_0 + h} \int_\Rd \left[ \alpha_t(x) \cdot \nabla \varphi(x) + \Delta \varphi(x) \right] d\mu_t(x)dt \right| \\
	&\le \norm{\nabla^2 \varphi}_\oo d^2_2(\mu_{t_0 + h}, \mu_{t_0}).
\end{align*}
The result then follows from part \eqref{prop:TimeLipschitzEstimate}, the assumption on $\alpha$, and the $L^1$-continuity in time of $\mu$ given by Lemma \ref{lem:L1cts}.

\end{proof}



The next ingredient we need is a differentiability property of the entropy along the flow of the Fokker-Planck equation. It will require preliminary facts on the continuity of the Fischer information.

\begin{lem}\label{lem:Fischerestimate}
	Assume $(b,c): [0,T] \times \RR^d \to \RR^d \times \RR$ are bounded, continuous, and continuously differentiable in the space variable with derivatives uniformly bounded on $[0,T] \times \RR^d$. Suppose that $u$ is a smooth, positive solution of
	\begin{equation}\label{genlin}
		\del_t u = \Delta u + b \cdot \nabla u + c u \quad \text{in } (t_0,T], \quad u(t_0,\cdot) = u_0,
	\end{equation}
	where $u_0 \in L^1_+(\Rd)$ and $\mcl I(u_0) < +\oo$. Then there exists a constant $C > 0$ depending only on $\norm{\nabla b}_\oo$ and $\norm{\nabla c}_\oo$ such that, for all $t \in (t_0,T]$,
	\[
		\mathcal{I}(u_t) \le e^{C(t-t_0)}\mathcal{I}(u_{t_0}) + (e^{C(t-t_0)} - 1) \sup_{s \in [t_0,t]} \norm{u(s,\cdot)}_{L^1}.
	\]
 
\end{lem}

\begin{proof}
	Define $F: \RR^d \times \RR_+ \to \RR$ by $F(v,u) = \frac{|v|^2}{u}$. Then $F$ is positive and convex, and so, for any smooth function $(v,u) : \RR^d \to \RR^d \times \RR_+$,
	\begin{align*}
		\Delta F(v,u) &= (\nabla_v F)(v,u) \cdot \Delta v + (\del_u F)(v,u) \Delta u 
		+ \sum_{i=1}^d 
		\left\langle
		D^2_{(v,u)} F(v,u)
		\begin{pmatrix}
			v_{x_i} \\
			u_{x_i}
		\end{pmatrix}
		,
		\begin{pmatrix}
			v_{x_i} \\
			u_{x_i}
		\end{pmatrix}
		\right\rangle
		\\
		&\ge (\nabla_v F)(v,u) \cdot \Delta v + (\del_u F)(v,u) \Delta u.
	\end{align*}
	Since $u$ is positive and smooth for $t > t_0$, the function $w := F(\nabla u, u)$ is well-defined, which therefore satisfies
	\begin{equation}\label{eq:forF}
	\begin{split}
		\del_t w& - \Delta w - b \cdot \nabla w- c w  \\
		&\le c\Big[ (\nabla_v F)(\nabla u, u) \cdot \nabla u + (\del_u F)(\nabla u, u) u - F(\nabla u, u) \Big] \\
		&+ (\nabla_v F)(\nabla u, u) \cdot \left[ (D b)\nabla u + (\nabla c) u \right]. 
	\end{split}
	\end{equation}

	Computing
	\[
		\del_u F(v,u) = - \frac{|v|^2}{u^2} \quad \text{and} \quad \nabla_v F(v,u) = \frac{2v}{u},
	\]
	we find that $v \cdot \nabla _v F + u \del_u F - F = 0$, and so \eqref{eq:forF} simplifies to
	\begin{align*}
		\del_t w - \Delta w - b \cdot \nabla w - c w
		&\le \frac{2Db \nabla u \cdot \nabla u}{u}+ 2 \nabla c \cdot \nabla u\\
		&\le 2 \norm{Db}_\oo w+ \norm{\nabla c}_\oo w + \norm{\nabla c}_\oo u,
	\end{align*}
	where in the last line we used Young's inequality. The result now follows easily upon integrating over $\Rd$, integrating by parts, and using Gr\"onwall's inequality. 
\end{proof}

Thanks to the above estimate we can show the following right-continuity property for the Fischer information along a solution to the Fokker-Planck equation.

\begin{lem} \label{lem:continuityFischer}
    Assume that $(\mu, \alpha)$ is a solution to the Fokker-Planck equation \eqref{FP} with $\mathcal{I}(\mu_0)<+\infty$ and $\alpha$ is smooth with bounded derivatives. Then $\mathcal{I}(\mu_t)<+\infty$ for all $t \geq t_0$ and $t\mapsto \mathcal{I}(\mu_t)$ is right-continuous.
\end{lem}


\begin{proof}
It suffices to check the right-continuity property at $t = t_0$. The right-continuity at any other point follows by a similar argument after the bound on $\mcl I(\mu_t)$ is established.

Rewriting the Fokker-Planck equation \eqref{FP} in non-divergence form,
the bound of Lemma \ref{lem:Fischerestimate} then implies, for all $t \geq t_0$ and some constant $C > 0$ depending on the $C^2$-norm of $\alpha$,
\begin{equation}\label{limsupI}
	\mathcal{I}(\mu_t) - \mathcal{I}(\mu_0) \leq \bigl( \mathcal{I}(\mu_0) +1 \bigr) \bigl( e^{C(t-t_0)}-1 \bigr).
\end{equation}
%
%
It follows that $(\nabla \sqrt{\mu_{t}})_{t > 0}$ is bounded in $L^2(\R^d)$. Moreover, by Lemma \ref{lem:L1cts},
	\[
		\int_{\R^d} |\sqrt{\mu_{t}} - \sqrt{\mu_0}|^2dx \le \int_{\RR^d} |\mu_{t} - \mu_0|dx \xrightarrow{t \to t_0^+} 0.
	\]
	This yields that, as $t \to t_0^+$, $\nabla \sqrt{\mu_{t}}$ converges weakly in $L^2(\R^d)$ toward $\nabla \sqrt{\mu_0}$. In particular, 
	\begin{equation}\label{liminfFischer}
		\liminf_{t \to t_0} \norm{\nabla \sqrt{\mu_{t}}}_{L^2(\Rd)} \ge \norm{\nabla \sqrt{\mu_0}}_{L^2(\R^d)},
	\end{equation}
	which, together with \eqref{limsupI}, gives the right-continuity of $\mcl I(\mu_t)$ at $t_0$. 

\end{proof}

\begin{prop}
    \label{prop:DiffEntropy} Assume that $(\mu, \alpha)$ is a solution to the Fokker-Planck equation \eqref{FP} with $\mathcal{I}(\mu_0)<+\infty$ and $\alpha$ smooth with bounded derivatives. Then
\[
	\mathcal{E}(\mu_t) = \mathcal{E}(\mu_0) + \int_{t_0}^t \int_{\R^d} \nabla \log\mu_s(x) \cdot \alpha_s(x) d\mu_s(x)ds - \int_{t_0}^t \mathcal{I}(\mu_s)ds, \quad \forall t \geq t_0.
\]
In particular, $t \mapsto \mathcal{E}(\mu_t)$ is differentiable at $t=t_0$ and 
\begin{equation}
\lim_{t \rightarrow t_0^+} \frac{\mathcal{E}(\mu_t) - \mathcal{E}(\mu_0)}{t-t_0} = \int_{\R^d} \nabla \log \mu_0 \cdot \alpha_{t_0}(x) d\mu_0(x) - \mathcal{I}(\mu_0).
\label{eq:differentialentropy}
\end{equation}
\end{prop}


\begin{proof}
	Define $e(s) = s\log s$ for $s \ge 0$. The smoothness of $\alpha$ and the properties of the Laplacian imply that, for $t > t_0$, $\mu_t$ is smooth and positive, so we can therefore justify the computations that lead to the equality, in the sense of distributions,
	\[
		\del_t e(\mu_t) + \div(\alpha_t e(\mu_t)) - \Delta e(\mu_t) = - (\div \alpha_t) \mu_t - \frac{|\nabla\mu_t|^2}{\mu_t}.
	\]
	Integrating in space and time (the former can be justified with cutoff functions of the form $\chi(x/R)$, then taking $R \to \oo$), we get
	\begin{equation}\label{Ediff:integral}
		\mcl E(\mu_t) - \mcl E(\mu_0) = - \int_{t_0}^t \int_{\RR^d} \left[ (\div \alpha_s) \mu_s + \frac{|\nabla\mu_s|^2}{\mu_s} \right]dx ds.
	\end{equation}
	where we have used that, by Lemma \ref{lem:Fischerestimate}, the last term is integrable. The derivative \eqref{eq:differentialentropy} then follows in view of the continuity of $t \to \mu_t$ in $\mathcal{P}_2(\R^d)$, the regularity of $\alpha$, and Lemma \ref{lem:continuityFischer}.
	\end{proof}

The final ingredient needed in the proof of the solution property is the following dynamic programming principle, whose proof, in this context, is straightforward since the dynamics are deterministic.

\begin{lem}[Dynamic Programming Principle] \label{L:DPP}
    For any $(t_0,\mu_0) \in [0,T) \times \mathcal{P}_2(\R^d)$ and $t_1 \in (t_0,T)$,
    $$U(t_0,\mu_0) = \inf_{(\mu,\alpha)} \left\{ \int_{t_0}^{t_1} \int_{\R^d} L \bigl(x, \alpha_t(x) \bigr)d\mu_t(x)dt + \int_{t_0}^{t_1} \mathcal{F}(\mu_t)dt + U \bigl(t_1, \mu_{t_1} \bigr) \right\},$$
where the infimum is taken over $(\mu, \alpha)$ with $\mu \in \mathcal{C} \bigl( [t_0,t_1], \mathcal{P}_2(\R^d) \bigl)$ and $\alpha \in L^2 \bigl([t_0,T] \times \R^d,\mu_t \otimes dt; \R^d \bigr)$ satisfying in the sense of distributions
$$ \partial_t \mu + \div( \alpha \mu) - \Delta \mu = 0,  \quad  \mbox{ in } (t_0,t_1) \times \R^d, \quad \quad \mu_{t_0} = \mu_0.$$
\end{lem}

\begin{rmk}\label{rmk:bndalpha}
	The $d_1$-Lipschitz continuity of $U$ given in Proposition \ref{prop: regularity of the value function} implies that the infimum in the dynamic programming principle can be restricted to controls $\alpha$ satisfying $\norm{\alpha}_{\infty} \leq C$, where $C$ is the same bound from that proposition.
 Moreover, in view of the $d_1$- and strong $L^1$-stability of the Fokker-Planck equation \eqref{FP} under perturbations of $\alpha$ implied by the presence of the Laplacian term, the infimum can be further restricted to smooth $\alpha$. This will be helpful in particular in the proof of the supersolution property.
\end{rmk}

\subsection{The subsolution property}

\begin{prop}
The value function $U$ defined in \eqref{def:Uexistence} is a viscosity sub-solution to the HJB equation \eqref{HJBExistence} in the sense of Definition \ref{D:soln} .
\end{prop}

\begin{proof}
Fix $\delta >0$ and take  $(t_0,\mu_0) \in [0,T) \times \mathcal{P}_2(\R^d)$ with $\mathcal{I}(\mu_0) <+\infty$ and $(r,p) \in \partial^+ \bigl( U - \delta \mathcal{E} \bigr) (t_0,\mu_0)$. Let $(\mu, \alpha)$ be admissible for the control problem with  $\alpha$ a smooth, bounded, time-independent control.  Using that $(r, p) \in \partial^+ \bigl(U - \delta \mathcal{E} \bigr) (t_0,\mu_0)$, Definition \ref{D:diff} gives, for all $t \in [t_0,T],$
\begin{align*}
	U(t,\mu_t) &\leq U(t_0,\mu_0) + r(t-t_0) + \int_{\R^d} p(x) \cdot \bigl( T_{\mu_0}^{\mu_t}(x) -x) d\mu_0(x) + o\bigl(d_2(\mu_t,\mu_0) +|t-t_0| \bigr) \\
 &+ \delta \bigl( \mathcal{E}(\mu_t) - \mathcal{E}(\mu_0) \bigr).
\end{align*}
Combined with the dynamic programming principle, Lemma \ref{L:DPP} and Propositions \ref{prop:FPresults} and \ref{prop:DiffEntropy}, we infer that
\begin{equation}\label{DPP18MAi2023}
	\begin{split}
		0 &\leq r(t-t_0) + \int_{\R^d} p(x) \cdot \bigl( T_{\mu_0}^{\mu_t}(x) -x) d\mu_0(x) + \int_{t_0}^t \int_{\R^d}L \bigl(x,\alpha(x)\bigr)d\mu_s(x)ds \\
		&+ \int_{t_0}^t \mathcal{F}(\mu_s)ds +\delta \int_{t_0}^t \int_{\R^d} \nabla \log \mu_s(x) \cdot\alpha(x)d\mu_s(x)ds - \delta \int_{t_0}^t \mathcal{I}(\mu_s)ds +  o\bigl(d_2(\mu_t,\mu_0) +|t-t_0| \bigr).
	\end{split}
\end{equation}
The definition of the tangent space (Definition \ref{D:tangent}) also implies that there exists a sequence $(\varphi_n)_{n \in \mathbb{N}} \subset \mathcal{C}_c^{\infty}(\R^d)$ such that, as $n \to \oo$, $\nabla \varphi_n$ converges to $p$ in $L^2(\mu_0)$. We then have, by the Cauchy-Schwarz inequality,
\[
	\int_{\R^d} p(x) \cdot \bigl( T_{\mu_0}^{\mu_t}(x) -x) d\mu_0(x) \leq \int_{\R^d} \nabla \varphi_n(x) \cdot \bigl( T_{\mu_0}^{\mu_t}(x) -x) d\mu_0(x) +  \|p - \nabla \varphi_n \|_{L^2(\mu_0)} d_2(\mu_t,\mu_0).
\]
By Proposition \ref{prop:FPresults}\eqref{prop:TimeLipschitzEstimate}, there exists $C > 0$ depending on $\norm{\alpha}_{L^\oo}$ and $\mcl I(\mu_0)$ such that, for all $t \in [t_0,T],$
\begin{equation}
d_2 \bigl(\mu_{t},\mu_0 \bigr) \leq C |t-t_0|.
\label{TimeLip18Mai2023}
\end{equation}
Inequality \eqref{DPP18MAi2023} then becomes, for some constant $C > 0$ independent of $t$ and $n$,
\begin{equation}\label{DPP218MAi2023}
	\begin{split}
	0 &\leq r(t-t_0) + \int_{\R^d} \nabla \varphi_n(x) \cdot \bigl( T_{\mu_0}^{\mu_t}(x) -x) d\mu_0(x) + \int_{t_0}^t \int_{\R^d}L \bigl(x,\alpha(x)\bigr)d\mu_s(x)ds \\
	&+ \int_{t_0}^t \mathcal{F}(\mu_s)ds  +\delta \int_{t_0}^t \int_{\R^d} \nabla \log \mu_s(x).\alpha(x)d\mu_s(x)ds - \delta \int_{t_0}^t \mathcal{I}(\mu_s)ds \\
 &+ C \|p - \nabla \varphi_n \|_{L^2(\mu_0)} |t-t_0| + o\bigl(|t-t_0| \bigr).
	\end{split}
\end{equation}
We then divide \eqref{DPP218MAi2023} by $t-t_0$ and send $t \rightarrow t_0^+$, which gives, upon appealing to Proposition \ref{prop:FPresults} and integrating by parts,
\begin{align*}
	0 &\leq r + \int_{\R^d} \bigl [ \nabla \varphi_n(x) \cdot \alpha(x) + L(x, \alpha(x)) + \Delta \varphi_n(x) \bigr ] d\mu_0(x) + C\|p - \nabla \varphi_n \|_{L^2(\mu_0)} \\
 &+ \delta \int_{\R^d} \nabla \log \mu_0(x) \cdot \alpha(x)d \mu_0(x) - \delta \mathcal{I}(\mu_0) \\
	&= r + \int_{\R^d}  \left[ \nabla \varphi_n(x) \cdot \alpha(x) + L(x, \alpha(x))  - \nabla \varphi_n(x) \cdot \frac{\nabla \mu_0(x)}{\mu_0(x)} \right] d\mu_0(x)+ C\|p - \nabla \varphi_n \|_{L^2(\mu_0)} \\
 &+ \delta \int_{\R^d} \nabla \log \mu_0(x) \cdot \alpha(x)d \mu_0(x) - \delta \mathcal{I}(\mu_0).
\end{align*}
Sending $n \to \oo$ and using the facts that $\nabla \varphi_n \rightarrow p$ in $L^2(\mu_0)$ and $\mathcal{I}(\mu_0)<+\infty$ gives
\begin{equation}\label{almostsub}
	0 \leq r + \int_{\R^d} \bigl [ \bigl(p(x) + \delta \nabla \log \mu_0(x) \bigr) \cdot \alpha(x)+ L(x, \alpha(x)) \bigr] d\mu_0(x) - \int_{\R^d} \nabla  \mu_0(x) \cdot p(x) dx - \delta \mathcal{I}(\mu_0).
\end{equation}
By Lemma \ref{lem:boundedp}, $\alpha^*(x) := -D_pH \bigl(x, p(x) + \delta \nabla \log \mu_0(x) \bigr)$ belongs to $L^\oo(\Rd)$. Approximating $\alpha^*$ in $L^2(\mu_0)$ by smooth and bounded controls and using the definition of $L$, we infer that
\begin{equation}
	0 \leq r - \int_{\R^d} H \bigl(x, \bigl(p(x) + \delta \nabla \log \mu_0(x) \bigr) \bigr)  d\mu_0(x) - \int_{\R^d} \nabla  \mu_0(x) \cdot p(x) dx - \delta \mathcal{I}(\mu_0).
\end{equation}
The local Lipschitz regularity of $H$ then implies that, for some constants $C,C'>0$ depending only on the Lipschitz constant of $U$ (see Lemma \ref{lem:boundedp} again),
\begin{align*}
	0 &\leq r - \int_{\R^d} H \bigl(x, p(x) \bigr)   d\mu_0(x) - \int_{\R^d} \nabla  \mu_0(x) \cdot p(x) dx  + C \delta \mathcal{I}(\mu_0)^{1/2} - \delta \mathcal{I}(\mu_0) \\
 & \leq r - \int_{\R^d} H \bigl(x, p(x)  \bigr)  d\mu_0(x) - \int_{\R^d} \nabla  \mu_0(x) \cdot p(x) dx  + C' \delta, 
\end{align*}
which gives the sub-solution property.

\end{proof}

\subsection{The supersolution property}

\begin{prop}
The value function $U$ defined in \eqref{def:Uexistence} is a viscosity super-solution to the HJB equation \eqref{HJBExistence} in the sense of Definition \ref{D:soln}.
\end{prop}

\begin{proof}
	Fix $\delta >0$ and take $(t_0,\mu_0)$ in $[0,T] \times \mathcal{P}_2(\R^d)$ such that $\mathcal{I}(\mu_0) <+\infty$ and $(r,p)\in \partial^- \bigl(U + \delta \mathcal{E} \bigr) (t_0,\mu_0)$. Fix $\epsilon > 0$ and let $(t_n)_{n \geq 1}$ be such that $t_n \rightarrow t_0^+$ as $n \rightarrow +\infty$. 
	By the dynamic programming principle Lemma \ref{L:DPP} we can find, for every $n \geq 1$, a couple $(\mu^n, \alpha^n)$ satisfying \eqref{FP} and such that
	\[
		U(t_0, \mu_0) \geq \int_{t_0}^{t_n} \int_{\R^d} L \bigl(x, \alpha_t^n(x) \bigr) d\mu_t^n(x)dt + \int_{t_0}^{t_n} \mathcal{F}(\mu_t^n)dt + U(t_n, \mu_{t_n}^n) - \epsilon(t_n-t_0).
	\]
	As noted in Remark \ref{rmk:bndalpha}), we can assume that, for all $n$, $\alpha^n$ is smooth with bounded derivatives, and, for some $C > 0$ independent of $n$ and $\epsilon$, $\norm{\alpha^n}_{\infty} \le C$. Therefore, there exists, by Proposition \ref{prop:FPresults}\eqref{prop:TimeLipschitzEstimate}, a constant $C>0$ independent of $n$ and $\epsilon$ such that
\begin{equation}
d_2(\mu^n_{t_n}, \mu_0) \leq C|t_n - t_0|, \quad \forall n \geq 1.
\label{d_2estimate5june2023}
\end{equation}
Therefore, using the definition of $\partial^{-} \bigl( U + \delta \mathcal{E} \bigr) (t_0, \mu_0)$, we find that, for all sufficiently large $n$,
\begin{align*}
	0 \geq r(t_n-t_0) &+ \int_{\R^d} p(x) \cdot \bigl(T_{\mu_0}^{\mu^n_{t_n}}(x)-x \bigr) d\mu_0(x) - \delta \bigl( \mathcal{E}(\mu_{t_n}^n) - \mathcal{E}(\mu_0) \bigr) \\
	&+ \int_{t_0}^{t_n} \int_{\R^d} L \bigl(x, \alpha_t^n(x) \bigr) d\mu_t^n(x)dt + \int_{t_0}^{t_n} \mathcal{F}(\mu_t^n)dt - 2\epsilon(t_n - t_0).
\end{align*}
In particular, for all $\varphi \in \mathcal{C}^{\infty}_c(\R^d)$, \eqref{d_2estimate5june2023} implies, by the Cauchy-Schwarz inequality, that
\begin{align*}
	0 &\geq r(t_n-t_0) + \int_{\R^d} \nabla \varphi(x) \cdot \bigl(T_{\mu_0}^{\mu^n_{t_n}}(x)-x \bigr) d\mu_0(x) - \delta \bigl( \mathcal{E}(\mu_{t_n}^n) - \mathcal{E}(\mu_0) \bigr)  \\
	&+ \int_{t_0}^{t_n} \int_{\R^d} L \bigl(x, \alpha_t^n(x) \bigr) d\mu_t^n(x)dt + \int_{t_0}^{t_n} \mathcal{F}(\mu_t^n)dt - \left( C \norm{ \nabla \varphi - p }_{L^2(\mu_0)} + 2\epsilon\right) (t_n-t_0).
\end{align*}
Using Proposition \ref{prop:DiffEntropy} and the regularity of $\alpha^n$ we have, using Cauchy-Schwarz inequality,
\begin{align*}
    \mathcal{E}(\mu_{t_n}^n) - \mathcal{E}(\mu_0) &= \int_{t_0}^{t_n} \int_{\R^d} \alpha^n_t(x) \cdot \nabla \log \mu_t^n(x) d\mu_t^n(x)dt - \int_{t_0}^{t_n} \mathcal{I}(\mu^n_t)dt \\
    & \leq \int_{t_0}^{t_n} \bigl( \norm{\alpha^n}_{\infty} \mathcal{I}(\mu^n_t)^{1/2} - \mathcal{I}(\mu^n_t) \bigr)dt \leq \frac{(t_n-t_0)}{4} \norm{\alpha^n}_{\infty}^2. 
\end{align*}
Appealing to Proposition \ref{prop:FPresults}\eqref{prop:IPP}, we obtain further that, for a constant $C > 0$ independent of $n$ and $\epsilon$,
\begin{align*}
	0 &\geq r(t_n-t_0) + \int_{t_0}^{t_n} \int_{\R^d} \bigl[ \alpha^n_t(x) \cdot \nabla \varphi(x) +\Delta \varphi(x) + L \bigl(x, \alpha_t^n(x) \bigr) \bigr] d\mu_t^n(x)dt \\
&+ \int_{t_0}^{t_n} \mathcal{F}(\mu_t^n)dt -\left( C \norm{ \nabla \varphi - p }_{L^2(\mu_0)} + \norm{\nabla^2 \varphi}_\oo(t_n - t_0) + C \delta + 2\epsilon \right)(t_n - t_0).
\end{align*}
Using the definition of $H_1$, this leads to
\begin{align*}
0 &\geq r(t_n-t_0) - \int_{t_0}^{t_n} \int_{\R^d} \bigl[ H_1 \bigl(x , \nabla \varphi(x) \bigr) -\Delta \varphi(x) \bigr] d\mu_t^n(x)dt + \int_{t_0}^{t_n} \mathcal{F}(\mu_t^n)dt \\
&-\left( C \norm{ \nabla \varphi - p }_{L^2(\mu_0)} + \norm{\nabla^2 \varphi}_\oo(t_n - t_0) + C \delta + 2\epsilon \right)(t_n - t_0).
\end{align*}
Dividing by $t_n-t_0$ and letting $n \rightarrow +\infty$, using the fact that $t \mapsto \mu_t^n$ is continuous in $L^1(\Rd) \cap \mcl P_1(\Rd)$, independently of $n$ and $\epsilon$, we arrive at
\begin{align*}
C \delta &\geq r - \int_{\R^d} \bigl[ H_1 \bigl(x, \nabla \varphi(x) \bigr)  \bigr]d\mu_0(x) + \mathcal{F}(\mu_0) -\int_{\R^d} \nabla \varphi(x) \cdot \nabla \mu_0(x)dx -C \norm{\nabla \varphi - p }_{L^2(\mu_0)} - 2 \epsilon.
\end{align*}
Finally, since $\varphi$ and $\epsilon$ were arbitrary, in view of Definition \ref{D:tangent}, we obtain
\[
	C \delta  \geq r - \int_{\R^d}  H_1 \bigl(x,p(x) \bigr) d\mu_0(x) + \mathcal{F}(\mu_0) -\int_{\R^d}  p(x) \cdot \nabla \mu_0(x) dx,
\]
which is the desired super-solution property.
\end{proof}

\bibliography{/Users/sam/Documents/DaudinSeeger2023ArxivSubmission/HJBrefs.bib}

\bibliographystyle{acm}

\end{document}